\documentclass[opre,nonblindrev]{informs3-bare}
\usepackage{amssymb,nicefrac,bm,upgreek,mathtools,verbatim,enumerate}
\usepackage[final]{hyperref}
\usepackage[mathscr]{eucal}
\RequirePackage[normalem]{ulem}
\usepackage{amsopn}

\usepackage{graphicx,subfigure,float}

\usepackage{natbib}
 \bibpunct[, ]{(}{)}{,}{a}{}{,}%
 \def\bibfont{\small}%

\TheoremsNumberedThrough     
\ECRepeatTheorems

\EquationsNumberedThrough    

\newcommand{\ttup}[1]{\textup{(}#1\textup{)}}
\newcommand{\stkout}[1]{\ifmmode\text{\sout{\ensuremath{#1}}}\else\sout{#1}\fi} 
\usepackage{color}
\definecolor{dmagenta}{rgb}{.4,.1,.5}
\definecolor{dblue}{rgb}{.0,.0,.5}
\definecolor{mblue}{rgb}{.0,.0,.8}
\definecolor{mred}{rgb}{.7,.0,.0}
\definecolor{dred}{rgb}{.6,.0,.0}
\definecolor{dred}{rgb}{.0,.0,.0}  
\definecolor{dgreen}{rgb}{.0,.4,.0}
\definecolor{Eeom}{rgb}{.0,.0,.5}

\newcommand{\rd}[1]{{\color{dred}#1}}

\hypersetup{
  colorlinks=true,
  citecolor=mblue,
  linkcolor=mblue,
  frenchlinks=false,
  pdfborder={0 0 0},
  naturalnames=false,
  hypertexnames=false,
  breaklinks}
\usepackage[capitalize,nameinlink]{cleveref}

\crefname{section}{Section}{Sections}
\crefname{subsection}{Section}{Sections}
\crefname{notation}{Notation}{Notations}
\crefname{condition}{Condition}{Conditions}
\crefname{hypothesis}{Hypothesis}{Conditions}
\crefname{assumption}{Assumption}{Assumptions}
\crefname{lemma}{Lemma}{Lemmas}
\crefname{fact}{Fact}{Facts}

\Crefname{figure}{Figure}{Figures}

\crefformat{equation}{\textup{#2(#1)#3}}
\crefrangeformat{equation}{\textup{#3(#1)#4--#5(#2)#6}}
\crefmultiformat{equation}{\textup{#2(#1)#3}}{ and \textup{#2(#1)#3}}
{, \textup{#2(#1)#3}}{, and \textup{#2(#1)#3}}
\crefrangemultiformat{equation}{\textup{#3(#1)#4--#5(#2)#6}}%
{ and \textup{#3(#1)#4--#5(#2)#6}}{, \textup{#3(#1)#4--#5(#2)#6}}%
{, and \textup{#3(#1)#4--#5(#2)#6}}

\Crefformat{equation}{#2Equation~\textup{(#1)}#3}
\Crefrangeformat{equation}{Equations~\textup{#3(#1)#4--#5(#2)#6}}
\Crefmultiformat{equation}{Equations~\textup{#2(#1)#3}}{ and \textup{#2(#1)#3}}
{, \textup{#2(#1)#3}}{, and \textup{#2(#1)#3}}
\Crefrangemultiformat{equation}{Equations~\textup{#3(#1)#4--#5(#2)#6}}%
{ and \textup{#3(#1)#4--#5(#2)#6}}{, \textup{#3(#1)#4--#5(#2)#6}}%
{, and \textup{#3(#1)#4--#5(#2)#6}}

\crefdefaultlabelformat{#2\textup{#1}#3}
\newcommand{\process}[1]{{\{#1_t\}_{t\ge0}}}

\newcommand{\dd}{\mathfrak{d}}
\newcommand{\Id}{\mathbb{I}}

\newcommand{\cB}{{\mathcal{B}}}  
\newcommand{\cE}{{\mathcal{E}}}  
\newcommand{\cG}{{\mathcal{G}}}  
\newcommand{\cI}{{\mathcal{I}}}  
\newcommand{\cJ}{{\mathcal{J}}}  
\newcommand{\Rm}{{\mathbb{R}^m}} 
\newcommand{\cK}{{\mathcal{K}}}  
\newcommand{\Lg}{\mathscr{A}}    

\newcommand{\sD}{{\mathscr{D}}}  

\newcommand{\bcZn}{{\Breve{\mathcal{Z}}^n}} 
\newcommand{\cZn}{{\mathcal{Z}^n}} 
\newcommand{\sS}{{\mathscr{S}}}  
\newcommand{\sX}{{\mathscr{X}}}  
\newcommand{\sZ}{{\mathscr{Z}}} 

\newcommand{\Lyap}{{\mathscr{V}}}  

\newcommand{\RR}{\mathbb{R}}
\newcommand{\NN}{\mathbb{N}}
\newcommand{\ZZ}{\mathbb{Z}}

\DeclareMathOperator{\Exp}{\mathbb{E}}

\newcommand{\D}{\mathrm{d}}

\newcommand{\Act}{{\mathbb{U}}}
\newcommand{\Usm}{\mathfrak{U}_{\mathrm{sm}}}

\newcommand{\Ind}{\bm{1}}   

\newcommand{\abs}[1]{\lvert#1\rvert}
\newcommand{\norm}[1]{\lVert#1\rVert}
\newcommand{\babs}[1]{\bigl\lvert#1\bigr\rvert}
\newcommand{\Babs}[1]{\Bigl\lvert#1\Bigr\rvert}
\newcommand{\babss}[1]{\biggl\lvert#1\biggr\rvert}
\newcommand{\bnorm}[1]{\bigl\lVert#1\bigr\rVert}

\newcommand{\transp}{^{\mathsf{T}}}
\newcommand{\df}{\coloneqq}

\DeclareMathOperator*{\diag}{diag}
\DeclareMathOperator*{\trace}{trace}

\newcommand{\order}{{\mathscr{O}}}
\newcommand{\sorder}{{\mathfrak{o}}}

\begin{document}

\TITLE{On System-Wide Safety Staffing 
of Large-Scale\\ Parallel Server Networks}
\MANUSCRIPTNO{}
\RUNTITLE{System-Wide Safety Staffing of Parallel Server Networks}
\RUNAUTHOR{H. Hmedi, A. Arapostathis,  and G. Pang}

\ARTICLEAUTHORS{%
\AUTHOR{Hassan Hmedi, Ari Arapostathis}
\AFF{Department of Electrical and Computer Engineering,
The University of Texas at Austin,\\ 2501 Speedway, EER 7.824,
Austin, TX~~78712,
\EMAIL{[hmedi,ari]@utexas.edu}}
\AUTHOR{Guodong Pang}
\AFF{The Harold and Inge Marcus Dept.
of Industrial and Manufacturing Eng.,
College of Engineering,\\
Pennsylvania State University,
University Park, PA~~16802,
\EMAIL{gup3@psu.edu}}
}

\ABSTRACT{
We introduce a ``system-wide safety staffing" (SWSS) 
parameter for multiclass multi-pool networks of any tree 
topology, Markovian or non-Markovian, in the Halfin--Whitt regime. 
This parameter can be regarded as the optimal 
reallocation of the capacity
fluctuations (positive or negative) of order $\sqrt{n}$ 
when each server pool employs 
a square-root staffing rule. 
We provide an explicit form of the SWSS as a function
of the system parameters, which is derived
using a graph theoretic approach based on Gaussian elimination.

For Markovian networks, we give an equivalent characterization of the SWSS 
parameter via 
the drift parameters of the limiting diffusion. 
We show that if the SWSS parameter is negative, the 
limiting diffusion
and the diffusion-scaled queueing processes are transient under any Markov control,  and cannot
have a stationary distribution when this parameter is zero.
If it is positive, we show that 
the diffusion-scaled queueing processes
are \emph{uniformly stabilizable}, that is,
there exists a scheduling policy under which the stationary distributions of the
controlled processes are tight
over the size of the network.
In addition, there exists a control under which the limiting controlled diffusion is
exponentially ergodic.
Thus we have identified a necessary and sufficient condition for the uniform
stabilizability of such networks in the Halfin-Whitt regime.  

We use a constant  control resulting from the leaf elimination algorithm to
stabilize the limiting controlled diffusion, while a family of Markov scheduling policies
which are easy to compute are used to stabilize the diffusion-scaled processes.
Finally, we show that under these controls the processes are exponentially ergodic and
the stationary distributions have exponential tails.}

\MSCCLASS{90B22, 60K25, 49L20, 90B36}


\maketitle

\section{Introduction}\label{S1} 

In recent years, parallel server networks have been a subject of intense study
due to their use in modeling a variety of systems including telecommunications, patient flows, 
service and data centers, etc.
The stability analysis of such systems is quite challenging because of their complexity. 
In this paper, we focus on studying the safety staffing and stability 
of such networks of any
tree topology in the Halfin--Whitt regime (or Quality--and--Efficiency--Driven (QED) regime)
in which the number of servers
and the arrival rates grow with the system scale while fixing the service
rates in a way that the system becomes critically loaded
\citep{H-W-81,w92,Borst-04}. 

When there is at least one class of jobs having a positive abandonment rate,
it is well known that 
there exists a scheduling policy under which the 
stationary distributions of
the controlled diffusion-scaled queueing processes are tight
(uniformly stable in the size of the network) \citep{AP16,AP18,AP19}. 
On the other hand, for networks with no abandonment, such results have only
been established for particular topologies. 
For the Markovian `V' network, it is shown in \citet{GS-12} and \citet{AHP18}
(the latter considers renewal arrivals, and the limiting diffusion) that, if
the server pool has $\sqrt{n}$ safety staffing, then
the stationary distributions of
the controlled diffusion-scaled queueing processes under
work-conserving stationary Markov scheduling policy are tight.
\rd{This is a very strong stability property since it is independent
of the system order or any particular work-conserving policy.
We say such networks are \emph{uniformly stable}.}
For the  `N' network,  \citet{Stolyar-15b} has shown that, with
$\sqrt{n}$ safety staffing in one server pool, the stationary distributions of
the controlled diffusion-scaled queueing processes are tight under 
a static priority scheduling policy. 
For a large class of Markovian networks,
which includes those with a single nonleaf server pool,
like the `N' and `M' models, and networks with class-dependent service rates,
a quantity referred to as $\sqrt{n}$ spare capacity is identified in \citet{HAP19},
and it is shown that when it is positive, the stationary distributions of
the controlled diffusion-scaled queueing processes are tight 
over the class of \emph{system-wide} work-conserving policies.
On the other hand, in \citet{Stol-Yud-13}, 
under a natural load balancing policy
referred to as ``Longest-Queue Freest-Server", 
it is shown that the stationary distributions of
the controlled diffusion-scaled queueing processes may not be tight 
for a network of arbitrary tree topology, but they are tight 
for the class of networks with pool-dependent service rates.

\rd{In Systems Theory, the existence of a control that renders a system stable, is
usually referred to as \emph{stabilizability}.
Adopting the same terminology, we say that a network is \emph{uniformly stabilizable}
if there exists some Markov scheduling policy under which the diffusion-scaled state process
is positive recurrent and the invariant distributions (over all
sufficiently large orders of the network) are tight.
The following question is then raised.
For parallel--server networks with an arbitrary tree topology and no abandonment,
is there a sharp criterion to determine if a network is uniformly stabilizable?
We are seeking a quantity
which if positive, the network is uniformly stabilizable, and if negative,
then the state process is transient under any Markov scheduling policy.
In this paper we provide an affirmative answer to the previous question through a
parameter called \emph{system-wide safety staffing} (SWSS), and which can
be easily computed from the system data.
Thus, the main result of the paper states that there exists a scheduling policy under which  
the stationary distributions of
the controlled  diffusion-scaled queueing processes are tight
if and only if the network has positive SWSS,
meaning that the SWSS parameter is positive.
As expected, for `V' and `N' networks, the SWSS parameter is the same as
the $\sqrt{n}$ safety staffing discussed in the preceding paragraph.}

For a better understanding of the SWSS parameter,
it is worth recalling the complete resource pooling (CRP)
condition and the demand and supply rates in the scales of order
$\order(n)$ and $\order(\sqrt{n})$  in the Halfin--Whitt regime.  
{\color{dred} The CRP condition requires that given a demand in the scale of order $\order(n)$,
each server pool has a suitable number of $\order(n)$ of servers 
so that there exists a unique allocation of the capacity in each server pool to meet
the demand of every class it can serve.
This of course determines the fluid limit.
More precisely, suppose that the arrival rates of the $n^{\mathrm th}$ system
are $\lambda^n_i = n\lambda_i$, where $i$ refers to the $i^\mathrm{th}$ class of jobs,
and the $\lambda_i$'s are fixed positive numbers.
Then the steady state allocations of servers in each pool are given by
the linear program \cref{E-LP} in \cref{S2}.
If the sever pools have an excess $\order(\sqrt{n})$ of servers from what is required
to meet this state allocation, then one can of course expect that
the system can be rendered stable by a suitable choice of a scheduling policy.
On the other hand, if some server pools are deficient, that is, they
are understaffed by an amount of $\order(\sqrt{n})$
servers, then the answer is not at all clear.
Thus, an important contribution of this paper, is that it quantifies
the `value' of a server in a given pool.
It can answer the question of whether moving a given number of servers from
one pool to another has a positive impact on system stability
(see \cref{R3}).
Moreover, the arrival rates in this paper also have a $\sqrt n$ component, that is,
$\lambda^n_i \approx n\lambda_i + \Hat{\lambda}_i \sqrt{n}$, and as a result
the SWSS parameter depends on the deviation of the arrival rates
from the nominal values $n\lambda_i$.
}

{\color{dred}
The SWSS parameter $\vartheta_p$ is obtained via the linear program \cref{E-LP'}, whereas
a similar program in \cref{E-LPn'} determines a parameter $\vartheta_p^n$ for the
$n^\mathrm{th}$ system (see \cref{D1}).
The asymptotic behavior of the system parameters in the Halfin--Whitt regime
(see \cref{E-Par01,E-Par02}) implies that $\vartheta_p^n$ tends to $\vartheta_p$
as $n\to\infty$.
It is asserted in \cref{T1} that if $\vartheta_p$ ($\vartheta_p^n$) are negative,
then the limiting diffusion ($n^\mathrm{th}$-system) are transient under
any Markov control. On the other hand, if $\vartheta_p>0$, then the limiting diffusion
is exponentially ergodic under some Markov control, and the $n^\mathrm{th}$-systems
are uniformly stabilizable for all large enough $n$.
Thus, the SWSS is an important and nontrivial extension of the familiar
square-root safety staffing parameter for single-class multi-server queues
\citep{H-W-81,w92}.}

A major contribution of this paper is a closed form expression for the SWSS
as a function of the system parameters.
Deriving this relies on solving the optimization
problem in \cref{E-LP'} via a simple Gaussian elimination of variables.
It is important to emphasize that the definition of
the SWSS and its functional form apply to
multiclass multi-pool networks of $G/G/N$ queues, regardless if they are 
Markovian or non-Markovian, since only the arrival and service rates play a role in this formulation. 

The results carry over to the limiting diffusion of the Markovian networks in an interesting manner.
We present in \cref{S4} a useful formula which allows
us to compute the SWSS as a function the drift parameters of the
limiting diffusion. 
{\color{dred} 
This relies on the explicit expression of the drift derived from using the iterative leaf elimination algorithm developed in \cite{AP16}, whose important properties are summarized in \cref{P2}. 
Moreover, we also provide an explicit matching expression (except an additional term indicating the violation of joint work conservation in the $n^{\rm th}$ system) for the infinitesimal drift of the 
diffusion-scaled Markovian queueing processes and some key properties of the main components in the expression in \cref{P1}. These properties of the drift expressions for both the diffusion-scaled processes and the limiting diffusions play a crucial rule for the  stability analysis. 
 }

{\color{dred} In \cref{S5} we show that the positivity of the SWSS is necessary for
stabilizability. In particular, we show in \cref{T5} that if $\vartheta_p<0$, then the
limiting diffusion process $\{X_t\}_{t\ge0}$ is transient under any Markov control,
and if $\vartheta_p=0$, then it cannot be positive recurrent.
Also, in \cref{T6}, we show that the exact analogous statement
(with the parameter $\vartheta^n_p$) applies to
the state process of the $n^{\mathrm{th}}$ system.
These results extend \citet[Propositions~3.1 and 3.2]{HAP19} to networks with
general tree topologies.
The proof of the above mentioned results relies on an important structural property of the drift
of the limiting diffusion stated in \cref{L2}
(see also \Cref{C2} for the drift of the diffusion-scaled state process).
}
 
In \cref{S6},
we exhibit a class of stabilizing controls for the diffusion-scaled queueing processes
and the limiting diffusion when the SWSS is positive. 
In order to accomplish this, we introduce an appropriate ``centering" for the
diffusion-scaled processes, which allows us to establish
Foster--Lyapunov equations. 
The stabilizing controls we use, consist of the family of balanced saturation policies
(BSPs) introduced in \citet{AP19}, where exponential ergodicity has been shown for
networks with at least one positive abandonment rates.
On the other hand, for the limiting diffusion, we use a constant control that relies
on the leaf elimination algorithm presented in \citet{AP16}, and show that it is also
stabilizing for the networks without abandonment.
We want to emphasize that the approach in \citet{AP16,AP19} does not apply to 
networks without abandonment.
We have focused on Markovian networks for the ease of exposition.
However, the stabilizability properties 
can be extended to networks with renewal arrivals and exponential service times 
using the methods in \citet{AHP18} (see \cref{R6,R9}). 

\smallskip

\noindent\emph{Organization of the paper.}
In the next subsection, we introduce the notation used in this paper.
In \cref{S2}, we describe the model, discuss the $\order(n)$ and  $\order(\sqrt{n})$
capacities, and introduce the SWSS parameter. 
In \cref{S3}, we present the calculation of the SWSS,
and provide the necessary and sufficient conditions on the fluctuations of order
$\order(\sqrt{n})$ to ensure that it is positive.
In \cref{S4}, we describe the system dynamics, introduce the re-centered
diffusion-scaled processes, and their diffusion limits.
We establish an equivalent characterization of the SWSS
in terms of the drift parameters and provide some examples.
In \cref{S5}, we establish the transience results both for the limiting
diffusion and diffusion-scaled processes in the case when the SWSS
is negative and show in addition that these processes cannot
be positive recurrent when this parameter is zero.
In \cref{S6.1}, we show that the limiting diffusion is exponentially ergodic under
a constant control.
In \cref{S6.2}, we prove that the BSPs are stabilizing, specifically, the
diffusion-scaled processes are exponentially ergodic under the BSPs. 

\subsection{Notation}\label{S1.1}
We use $\mathbb{R}^m$ (and $\mathbb{R}^m_+$), $m\ge 1$,
to denote real-valued $m$-dimensional (nonnegative) vectors, and write $\RR$
for the real line.
The transpose of a vector $z\in\Rm$ is denoted by $z\transp$.
Throughout the paper, $e\in\Rm$ stands for the vector whose elements are
equal to $1$, that is, $e=(1,\dotsc,1)\transp$, and $e_i\in\Rm$ denotes
the vector whose elements are all $0$ except for the $i^{\mathrm{th}}$ element
which is equal to $1$.
For a set $A\subseteq\RR^m$, we use $A^c$, and $\Ind_{A}$
to denote the complement, and the indicator function of $A$, respectively.
The Euclidean norm on $\RR^m$ is denoted by $\abs{\,\cdot\,}$,
and $\langle \cdot\,,\,\cdot\rangle$ stands for the inner product.
For a finite signed measure $\nu$ on $\Rm$,
and a Borel measurable $f\colon\Rm\to[1,\infty)$,
the $f$-norm of $\nu$ is defined by
\begin{equation*}
\norm{\nu}_f \,\df\, \sup_{\substack{g\in\cB(\Rm), \; \abs{g}\le f}}\;
\babss{\int_{\Rm} g(x)\,\nu(\D{x})}\,,
\end{equation*}
where $\cB(\Rm)$ denotes the class of Borel measurable functions on $\Rm$.

\section{Model Description and Summary of the Results} \label{S2}
We study multiclass multi-pool Markovian networks with $I$
classes of customers and $J$ server pools, and let
$\cI = \{1,\dotsc,I\}$ and $\cJ = \{1, \dotsc, J\}$. 
Customers of each class form their own queue, are served in the
first-come-first-served (FCFS) service discipline, and
do not abandon/renege while waiting in queue.
The buffers of all classes are assumed to have infinite capacity.  
We assume that the customer arrival and service processes of
all classes are mutually independent. 
We let $\cJ(i) \subset \cJ$, denote the subset of server pools that can serve class
$i$ customers, and $\cI(j) \subset \cI$ the subset
of customer classes that can be served by server pool $j$.
We form a bipartite graph $\cG = (\cI\cup \cJ, \cE)$ with a set
of edges defined by $\cE = \{(i,j)\in\cI\times\cJ\colon j\in\cJ(i)\}$, and
use the notation $i \sim j$, if $(i,j)\in\cE$,
and $i \nsim j$, otherwise.
We assume that the graph $\cG$ is a tree.

We consider a sequence of such network systems with the associated variables,
parameters and processes indexed by $n$. 
We study these networks in the Halfin--Whitt regime
(or the Quality-and-Efficiency-Driven (QED) regime), where the arrival
rate of each class and the number of servers in each pool grow
large as $n \to \infty$ in such a manner that the system becomes
critically loaded.
Note that the model description and the asymptotic regime apply to both
Markovian and non-Markovian networks. 

Let $\lambda_i^n$ and $\mu^n_{ij}$
be positive real numbers denoting the arrival rate of class-$i$
and the service rate of class-$i$ at pool $j$ if $i\sim j$
in the $n^{\rm th}$ system, respectively. 
Also $N^n_j$ is a positive integer denoting the number of servers in pool $j$.
The standard assumption concerning these parameters in the Halfin--Whitt regime
is that the following limits exist as $n\to\infty$:
\begin{equation}\label{E-Par01}
\frac{\lambda^{n}_{i}}{n} \,\to\, \lambda_i \,>\, 0\,,\quad
\frac{N^{n}_{j}}{n} \,\to\, \nu_{j} \,>\, 0\,,\quad
\mu_{ij}^n \,\to\, \mu_{ij} \,>\, 0\,,
\end{equation}
\begin{equation}\label{E-Par02}
 \frac{\lambda^{n}_{i} - n \lambda_i}{\sqrt n} \,\to\, \Hat\lambda_i\in\RR\,,
\quad
{\sqrt n}\,(\mu^{n}_{ij} - \mu_{ij}) \,\to\, \Hat\mu_{ij}\in\RR\,,\quad
\text{and}\quad {\sqrt n}\,(n^{-1} N^{n}_{j} -  \nu_{j})
\,\to\, \Hat\nu_{j}\in\RR\,.
\end{equation}

Let
\begin{equation*}
\RR^{\cG}_+ \,\df\, \bigl\{\xi=[\xi_{ij}]\in\RR^{I \times J}_+\colon
\xi_{ij}=0~~\text{for~}i\nsim j\bigr\}\,,
\end{equation*}
and analogously define $\RR^{\cG}$, $\ZZ_+^{\cG}$, and $\ZZ^{\cG}$.
We assume that the \emph{complete resource pooling} (CRP) condition is satisfied
(see \cite{Williams-00, Atar-05b}),  that is, 
the linear program (LP) given by
\begin{equation}\label{E-LP}
\begin{aligned}
&\text{minimize} \quad  \max_{j \in \cJ}\;\sum_{i \in \cI} \xi_{ij}\quad
\text{over\ \ } [\xi_{ij}] \in\RR^{\cG}_+\,,\\
\quad
&\text{subject to} \quad  \sum_{j \in \cJ} \mu_{ij} \nu_j \xi_{ij}
\,=\,  \lambda_i\quad \forall\, i \in \cI\,,
\end{aligned}\tag{\sf{LP}}
\end{equation}
has a unique solution
$\xi^*=[\xi^*_{ij}]\in\RR^{\cG}_+$  satisfying 
\begin{equation} \label{E-critfl}
\sum_{i \in \cI} \xi^*_{ij} \,=\, 1, \quad \forall j \in \cJ \,,
\quad\text{and}\quad \xi^*_{ij}\,>\,0\quad \text{for all~} i \sim j\,.
\end{equation}

We define  $x^* = (x_{i}^{*})_{i \in \cI}\in\RR^I_+\,,$ and
$z^* = [z_{ij}^*]\in\RR^{\cG}_+$ by
\begin{equation} \label{E-statfl}
x_{i}^{*} \,\df\,\sum_{j \in \cJ} \xi^*_{ij} \nu_j\,,
\qquad z_{ij}^* \,\df\, \xi_{ij}^* \nu_j \,. 
\end{equation}
The variable  $x^*_i$ can be interpreted as the steady-state number
of customers in class $i$, and the variable $z^*_{ij}$ as the steady-state number
of customers in each class $i$ receiving service in pool $j$, in the fluid scale.
Note that the steady-state queue lengths are all zero in the fluid scale.
The quantity $\xi_{ij}^*$  can be interpreted as the steady-state fraction of
service allocation of pool~$j$ to class-$i$ jobs in the fluid scale.
It is evident that \cref{E-critfl,E-statfl} imply that
$\sum_{i\in\cI} x^*_i=  \sum_{j\in\cJ}   \nu_j$.
For more details on this model, we refer the reader to
\citet{Atar-05a,Atar-05b} and \citet{AP16,AP19}.

\subsection{The System-Wide Safety Staffing Parameter}

Let $\{N_{ij}^n\in\NN\,,\;(i,j)\in\cE\,,n\in\NN\}$ be a sequence which satisfies
\begin{equation}\label{E-Nij}
\lfloor \xi_{ij}^* N_j^n \rfloor \,\le\, N_{ij}^n
\,\le\, \lceil \xi_{ij}^* N_j^n \rceil\,,\qquad\text{and\ \ }
\sum_{i\in \cI(j)}N_{ij}^n \,=\, N_j^n\,.
\end{equation}
By \cref{E-Par01,E-Par02}, we can write
\begin{equation}\label{E-apprA}
N^n_{ij} \, = \, z^*_{ij} n +  \xi^*_{ij} \Hat\nu_j \sqrt{n}
+ \sorder\bigl(\sqrt{n}\bigr) \,,
\quad (i,j) \in \cE\,,
\end{equation}
where we use the definition in \cref{E-statfl}.
Similarly, we have
\begin{equation}\label{E-apprB}
\lambda^n_i \, = \, \lambda_i n +\Hat\lambda_i \sqrt n
+ \sorder\bigl(\sqrt{n}\bigr)\,,\quad\text{and\ \ }
\mu^n_{ij} \, = \, \mu_{ij} + \frac{\Hat\mu_{ij}}{\sqrt n}
+ \sorder\biggl(\frac{1}{\sqrt n}\biggr)\,.
\end{equation}
By combining \cref{E-apprA,E-apprB} and the constraint in the (LP), we obtain 
\begin{equation}\label{E-spareA}
\sum_{j \in \cJ(i)} \mu^n_{ij} N^n_{ij} - \lambda^n_i \, = \,
- \Hat\lambda_i\sqrt n +
\sqrt{n} \sum_{j \in \cJ(i)}
\bigl(\mu_{ij} \xi^*_{ij}\Hat\nu_j  + \Hat\mu_{ij} z^*_{ij}\bigr)
+ \sorder\bigl(\sqrt{n}\bigr)
\qquad \forall\, i \in \cI\,. 
\end{equation}
Thus, for class $i$ customers,
the total steady-state servers allocated from all pools 
may be deficient, or have a surplus, of order $\order(\sqrt{n})$.

Recall that in the single class, single pool case (with $N$ servers)
the safety staffing parameter $\vartheta$
is given by
\begin{equation}\label{E-safe}
N \,=\, \nicefrac{\lambda}{\mu} + \vartheta \sqrt{\nicefrac{\lambda}{\mu}}\,.
\end{equation}
\rd{Let $\varDelta_I$ denote the set of probability vectors in $\RR^I$,
and $p = (p_1,\dotsc,p_I)$ be a positive vector in $\varDelta_I$.}
Mimicking \cref{E-safe}, to extend the definition
of the safety staffing parameter to the multiclass, multi-pool case, we seek
an alternate set of allocations $\{\widetilde{N}^n_{ij}\in\NN\colon i\sim j\}$
satisfying
\begin{equation}\label{E-spareB}
\begin{aligned}
\sum_{i\in \cI(j)}\widetilde{N}_{ij}^n
&\,=\, N_j^n\quad\forall\, j\in\cJ\,,\quad\text{and} \\
\sum_{j \in \cJ(i)} \mu^n_{ij} \widetilde{N}^n_{ij} - \lambda^n_i
\, & = \,  \vartheta_p\,p_i \sqrt{n}  + \sorder\bigl(\sqrt{n}\bigr)
\quad \forall\, i \in \cI\,,
\end{aligned}
\end{equation}
for some constant $\vartheta_p$.
If \cref{E-spareB} holds for some $\vartheta_p>0$ and a positive vector $p\in\varDelta_I$,
then as we show in
\cref{T8}, the system is uniformly stabilizable in the sense of the definition in \cref{S1}.

It is clear by \cref{E-spareA,E-spareB} and the complete resource pooling hypothesis,
that $\abs{N^n_{ij}-\widetilde{N}_{ij}^n}= \order(\sqrt n)$.
Thus $\widetilde{N}_{ij}^n$ has the form
\begin{equation}\label{E-spareC}
\widetilde{N}_{ij}^n \,=\, z^*_{ij} n + \Tilde{\kappa}_{ij} \sqrt n
+ \sorder\bigl(\sqrt{n}\bigr)
\end{equation}
for some $\Tilde{\kappa}=[\Tilde{\kappa}_{ij}]\in\RR^{\cG}$.
By \cref{E-apprA,E-spareB,E-spareC}, we have
$\sum_{i\in\cI} \Tilde{\kappa}_{ij} = \Hat\nu_j$.
It also follows from \cref{E-spareA} that such a collection
$\widetilde{N}^n_{ij}$ satisfying \cref{E-spareB} with $\vartheta_p>0$ can be found if
and only if
the linear program
\cref{E-LP'} in \cref{D1} below has a positive solution $\vartheta_p$.

\begin{definition}\label{D1}
We use a positive $p\in\varDelta_I$ as a free parameter.
Abusing the notation, let $\vartheta_p$ and $\kappa=[\kappa_{ij}]\in\RR^{\cG}$
be the unique  solution to  the linear program:
\begin{equation}\label{E-LP'}
\begin{aligned}
\text{maximize}&\quad \vartheta_p \\
\text{subject to}&\quad
\Hat\lambda_i\,\le\,
\sum_{j \in \cJ(i)} \mu_{ij}\kappa_{ij}-\vartheta_p\,p_i 
\qquad\forall\,i\in\cI\,,\\
& \sum_{i\in\cI(j)} \kappa_{ij} \,=\, \theta_j\,\df\,\Hat\nu_j +
 \sum_{i\in\cI(j)}\frac{\Hat\mu_{ij}}{\mu_{ij}} z^*_{ij}
 \qquad\forall j\in\cJ\,.
 \end{aligned}\tag{\sf{LP$^\prime$}}
\end{equation}
We refer to $\vartheta_p$ as the SWSS parameter, or simply as the SWSS.

We also define $\vartheta^n_p$
and $\kappa^n=[\kappa^n_{ij}]\in\RR^{\cG}$ as the unique solution to 
\begin{equation}\label{E-LPn'}
\begin{aligned}
\text{maximize}&\quad \vartheta^n_p \\
\text{subject to}&\quad
\Hat\lambda^n_i\,\le\,
\sum_{j \in \cJ(i)} \mu^n_{ij}\kappa^n_{ij}-\vartheta^n_p\, p_i 
\qquad\forall\,i\in\cI\,,\\
& \sum_{i\in\cI(j)} \kappa^n_{ij} \,=\, \theta^n_j\,\df\,\Hat\nu^n_j +
\sum_{i\in\cI(j)}\frac{\Hat\mu^n_{ij}}{\mu^n_{ij}} z^*_{ij}
 \qquad\forall j\in\cJ\,,
 \end{aligned}\tag{$\mathsf{LP}_n^\prime$}
\end{equation}
with
\begin{equation*}
    \Hat\lambda_i^n\,\df\, \frac{\lambda_i^n - n\lambda_i}{\sqrt{n}}\,,
    \quad \Hat\nu_{j}^n\,\df\,{\sqrt n}\,(n^{-1} N^{n}_{j} -  \nu_{j})\,, \quad 
    \Hat\mu_{ij}^n \,\df\, {\sqrt n}\,(\mu^{n}_{ij} - \mu_{ij})\,.
\end{equation*}
\end{definition}

Note that the CRP condition consists
of solving the first-order optimization problem \cref{E-LP}
(the quantities of order $n$,
matching supply and demand in the fluid scale), while \cref{E-LP'}
can be regarded as a second-order optimization problem
(the quantities of order $\order(\sqrt{n})$
involved in the `reallocation' of staffing).
Note also that $\Hat\lambda_i^n$, $\Hat\mu_{ij}^n$,
and $\Hat\nu_j^n$ converge to $\Hat\lambda_i$,
$\Hat\mu_{ij}$, and $\Hat\nu_j$ respectively as $n\to\infty$ by \cref{E-Par02}.

\rd{\begin{remark}\label{R1}
We note here that the sign of \rd{$\vartheta_p$} does not depend on
the positive vector $p$ chosen.
The proof of this fact is clear from the statement of \cref{T2}. In addition, we note that the choice of the vector $p$ plays a crucial role in the proof of stability of the limiting diffusion.
In particular, in the proof of \cref{T7} we have to select
$p$ such that $\langle p, S e_I\rangle>0$, where $S$ is a positive definite matrix.
This is the primary reason behind the introduction of the vector $p$.
But note also the identities in \cref{T4,R7}.
\end{remark}
}

{\color{dred}
\begin{remark}
The uniqueness of the solutions to \cref{E-LP',E-LPn'}
follows from the tree structure.
In fact, if we replace the inequality in the constraint of \cref{E-LP'}
with equality,
then we obtain $I+J$ independent equations
in the variables $[\kappa_{ij}]\in\RR^{\cG}$ and $\vartheta_p$,
and the same applies to \cref{E-LPn'}.
Thus these linear programs are equivalent to a system of linear equations.
The reason that we write them in this form is because as it follows from
the proof of \cref{T8}, that any feasible solution of \cref{E-LP'} with $\vartheta_p>0$ can
be used to synthesize a stabilizing scheduling policy.
\end{remark}
}

\subsection{Summary of the Results} 
In \cref{S3} we solve for $\vartheta_p$ as a function of the system parameters.
There is another significant result which is established in \cref{S4}.
As shown in \citet{AP16}, the drift of the limiting diffusion of
Markovian parallel server networks has the form
\begin{equation*}
b(x, u) \,=\, h - B_1 (x -  \langle e, x \rangle^{+} u^c)
+ \langle e, x \rangle^{-} B_2 u^s\,,
\end{equation*}
where $B_1$ and $B_2$ are in $\RR^{ I\times I}$ and $\RR^{I\times J}$,
respectively (see \cref{P2}).
Also, the vector $h=(h_i)_{i\in\cI}$ is given by (compare with \cref{E-spareA})
\begin{equation}\label{E-h}
h_i \,\df\, \Hat\lambda_i -
 \sum_{j \in \cJ(i)}
\bigl(\mu_{ij} \xi^*_{ij}\Hat\nu_j  + \Hat\mu_{ij} z^*_{ij}\bigr)\,,
\qquad i \in \cI\,.
\end{equation}
As shown in \citet{HAP19} the quantity
$\varrho \df - \langle e, B_1^{-1} h\rangle$ characterizes the uniform
stability of multiclass multi-pool networks that have a single
non-leaf server node (such as the `M' network) or those with class-dependent
service rates.
For this class of networks, it is shown that the system has an invariant probability
distribution under any stationary Markov control (i.e., uniformly stable)
if and only if $\varrho>0$.
The parameter $\varrho$ is referred to as `spare capacity' in that paper. 
(It is worth mentioning that this spare capacity is also used for the stability
of diffusions with jumps arising from many-server queues with abandonment in \citet{AHPS19,APS19,APS20}).
We show in \cref{S4} that for any multiclass multi-pool network with the
above diffusion limit, 
it holds that
$\varrho = \langle e,B_1^{-1} p \rangle\, \vartheta_p$.
Then, we show that $\vartheta_p>0$ is
a \emph{necessary and sufficient condition} for a multiclass multi-pool network
as described above to be stabilizable.
This also applies to the diffusion-scaled processes.
In fact we show that there exists a suitable scheduling
policy that renders the processes exponentially ergodic.
This result is summarized in the following theorem,
whose proof follows from \cref{T5,T6,T8,T7}.

\begin{theorem}\label{T1}
The following hold:
\begin{itemize}
\item[\ttup a]
If $\vartheta_p>0$, then the diffusion-scaled processes
and the limiting diffusion
are stabilizable.  Moreover, there exists a family of Markov scheduling policies,
under which the diffusion-scaled processes are exponentially ergodic and their
stationary distributions are tight and have exponential tails for all sufficiently
large system orders.
The same is true for the limiting diffusion under some stationary Markov control.
\item[\ttup b]
If $\vartheta_p<0$ ($\vartheta_p=0$), then the
limiting diffusion is transient (cannot have an invariant
probability measure) under any stationary Markov control.
The same applies to the state process of the $n^{\mathrm th}$ system
with respect to $\vartheta^n_p$ for all $n>0$.
\end{itemize}
\end{theorem}

\section{Computing the SWSS Parameter}\label{S3}

\Cref{T2} below, provides an explicit solution to
\cref{E-LP'}.
For this, we need some additional notation.
Let $(i,j)\in\cI\times\cJ$.
With $(i_1, j_1, i_2, j_2,\dotsc, i_m,j_m)$ denoting
the unique path of minimum length connecting $i\equiv i_1$ to $j\equiv j_m$ in $\cG$, 
we define the ``gain'' $\dd(i,j)$ by
\begin{equation*}
\dd(i,j) \,\df\, \mu_{i_1j_1}
\prod_{k=1}^{m-1} \frac{\mu_{i_{k+1}j_{k+1}}}{\mu_{i_{k+1}j_k}}\,.
\end{equation*}
Similarly, we define the gain $\dd(i,i')$
between any pair $i, i'\in \cI$, $i\ne i'$, by
\begin{equation}\label{E-dd}
\dd(i,i') \,\df\,
\prod_{k=1}^{m-1} \frac{\mu_{i_{k}j_{k}}}{\mu_{i_{k+1}j_k}}\,,
\end{equation}
where the product in \cref{E-dd} is evaluated over the
analogous path  $(i_1, j_1, i_2, j_2,\dotsc, i_m)$
connecting $i\equiv i_1$ to $i'\equiv i_m$ in $\cG$,
and we let $\dd(i,i)\df1$ for $i\in\cI$.

\begin{theorem}\label{T2}
The solution $\vartheta_p$ to \cref{E-LP'} is given by 
\begin{equation}\label{ET2A}
\vartheta_p  \,=\, \frac{\sum_{j\in\cJ} \dd(i,j) \theta_j
- \sum_{\ell\in\cI} \dd(i,\ell) \Hat\lambda_{\ell}}
{\sum_{\ell\in\cI}\dd(i,\ell)p_\ell}
\qquad \forall\,i\in\cI\,.
\end{equation}  
\end{theorem}


\proof{Proof.}
Consider a network graph $\cG$ as described in \cref{S2},
and a set of parameters $\Theta = \{\theta_j\colon j\in\cJ\}$.
Suppose there exist $\vartheta_p\in\RR$,
and a collection $\mathfrak{K} = \{\kappa_{ij} \colon (i,j)\in\cE\}$
solving
\begin{equation}\label{PT2A}
\sum_{j\in\cJ(i)} \mu_{ij}\kappa_{ij} \,=\,\Hat\lambda_i + \vartheta_p\, p_i 
\quad \forall\,i\in\cI\,,
\quad\text{and}\quad
\sum_{i\in \cI(j)} \kappa_{ij}\,=\, \theta_j\quad \forall\,j\in\cJ\,.
\end{equation}
We use $\mathfrak{K}(\cG,\Theta)$ and
$\kappa_{ij}(\cG,\Theta)$
to indicate explicitly the dependence of the solution on
the graph and the parameters $\Theta$.
The parameters $p$, $\vartheta_p$, and $\Hat\lambda=(\Hat\lambda_i)_{i\in\cI}$
are held fixed throughout
the proof.  

Let $\cI_{\mathsf{leaf}}$ denote the customer classes in $\cI$ which are
leaves of the graph. 
Consider the subgraph $\cG^0 = (\cI^0 \cup \cJ^0,\cE^0)$, with
$\cI^0=\cI\setminus\cI_{\mathsf{leaf}}$, $\cJ^0=\cJ$, and
$\cE^0 = \{(i,j)\in\cE \colon i\in\cI^0,j\in\cJ^0\}$.
{\color{dred}
Let $\Theta^0=\{\theta_j^0\colon j\in\cJ\}$, where
\begin{equation}\label{PT2B}
\theta_j^0 \,=\,
\theta_j - \sum_{i\in\cI_{\mathsf{leaf}}\cap\cI(j)} 
\mu_{ij}^{-1}\bigl(\Hat\lambda_i+\vartheta_p\,p_i\bigr)\,,\quad j\in\cJ\,.
\end{equation}
We claim that \cref{PT2A} has a solution for $(\cG,\Theta)$ if and only
if it is solvable for $(\cG^0,\Theta^0)$, and that
\begin{equation}\label{PT2Bb}
\kappa_{ij}(\cG,\Theta)\,=\, \kappa_{ij}(\cG^0,\Theta^0)\qquad
\forall\,(i,j)\in\cE^0\,.
\end{equation}
To prove the claim, let $\kappa_{ij}= \kappa_{ij}(\cG,\Theta)$ be solution for
$(\cG,\Theta)$.
It is clear from \cref{PT2A} that
$\kappa_{ij}= \mu_{ij}^{-1}\bigl(\Hat\lambda_i+\vartheta_p\,p_i\bigr)$
for $i\in\cI_{\mathsf{leaf}}$.
Since $\cI(j)$ is the disjoint union of $\cI^0(j)$ and $\cI_{\mathsf{leaf}}\cap\cI(j)$,
we write the second equation in \cref{PT2A} as
\begin{equation*}
\sum_{i\in\cI^0(j)}\kappa_{ij} \,=\,
\theta_j - \sum_{i\in\cI_{\mathsf{leaf}}\cap\cI(j)} \kappa_{ij}
\,=\, \theta_j - \sum_{i\in\cI_{\mathsf{leaf}}\cap\cI(j)}
\mu_{ij}^{-1}\bigl(\Hat\lambda_i+\vartheta_p\,p_i\bigr)\,=\, \theta^0_j \,.
\end{equation*}
It is also clear from the definitions that
$\sum_{j\in\cJ^0(i)} \mu_{ij}\kappa_{ij} \,=\,\vartheta_p\, p_i +\Hat\lambda_i$
for all $i\in\cI^0$.
Thus we obtain a solution for $(\cG^0,\Theta^0)$ as claimed.
Conversely if we start from a solution $\kappa_{ij}(\cG^0,\Theta^0)$, and augment
this by
defining $\kappa_{ij}= \mu_{ij}^{-1}\bigl(\vartheta_p\, p_i + \Hat\lambda_i\bigr)$
for $\cI_{\mathsf{leaf}}$,
we obtain a solution for $(\cG,\Theta)$.

}


Continuing, we claim that for any network graph $\cG$ which contains no customer leaves
and $\cI$ is not a singleton
there exists some $i\in\cI$ such that $\cJ(i)$ contains exactly one non-leaf element.
If the claim were not true, then removing all server leaves would result in a graph that
has no leaves, which is impossible since the resulting graph has to be a nontrivial tree.

Suppose then that $\cI^0$ is not a singleton, otherwise we are at the last step
of the construction which we described next.
Let $\imath^{}_1\in\cI^0$ be such that
exactly one member of $\cJ(\imath^{}_1)$, denoted as $\jmath^{}_1$,  is a
non-leaf in $\cG^0$.
Define
\begin{equation}\label{PT2C}
\theta_j^1 \,=\, \begin{cases} \theta_{\jmath^{}_1}^0
- \mu_{\imath^{}_1\jmath^{}_1}^{-1}
\Bigl(\Hat\lambda_{\imath^{}_1}+
\vartheta_p\, p_{\imath_1}
- \sum_{k\in\cJ(\imath^{}_1)\setminus\{\jmath^{}_1\}}
\mu_{\imath^{}_1 k}\,\theta_k^0\Bigr) & \text{for~} j=\jmath^{}_1\,,\\
\theta_j^0 & \text{for~} j \ne \jmath^{}_1\,.
\end{cases}
\end{equation}
Let $\cG^1=(\cI^1 \cup \cJ^1,\cE^1)$
denote the subgraph of $\cG^0$ which arises if we remove all
the edges containing $\imath^{}_1$ from $\cE^0$,
and define $\Theta^1 \df \{\theta_j^1\colon j\in\cJ^1\}$.
By \cref{PT2A}, we have
\begin{equation}\label{PT2D}
\kappa_{\imath^{}_1\jmath^{}_1}(\cG^0,\Theta^0)
\,=\, \mu_{\imath^{}_1\jmath^{}_1}^{-1}
\Biggl(\Hat\lambda_{\imath^{}_1} +
\vartheta_p\, p_{\imath_1} -
\sum_{k\in\cJ(\imath^{}_1)\setminus\{\jmath^{}_1\}}
\mu_{\imath^{}_1 k}\,\theta_k^0\Biggr)\,.
\end{equation}
It is clear then by \cref{PT2C,PT2D} that
\cref{PT2A} has a solution for $(\cG^0,\Theta^0)$ if and only
if it is solvable for $(\cG^1,\Theta^1)$.

Iterating the procedure in the preceding paragraph we obtain
a decreasing sequence of subgraphs
$\cG^\ell=(\cI^\ell \cup \cJ^\ell,\cE^\ell)$
for $\ell=1,\dotsc,m\df\abs{\cI^0}-1$, such that
$\cI^{m}$ is a singleton,
together with a sequence of parameter sets
$\Theta^\ell \df \{\theta_j^\ell\colon j\in\cJ^\ell\}$
and pairs $(\imath_\ell,\jmath_\ell)\in\cE^{\ell}$,
satisfying
\begin{equation}\label{PT2E}
\theta_j^{\ell} \,=\,
\begin{cases} \theta_{\jmath_\ell}^{\ell-1}
- \mu_{\imath_\ell\jmath_\ell}^{-1}
\Bigl(\Hat\lambda_{\imath_\ell} + \vartheta_p\, p_{\imath_\ell}
- \sum_{k\in\cJ(\imath_\ell)\setminus\{\jmath_\ell\}}
\mu_{\imath_\ell k}\,\theta_k^{\ell-1}\Bigr)
& \text{if~} j=\jmath_\ell\,,\\
\theta_j^{\ell-1}
& \text{if~} j \ne \jmath_\ell\,,
\end{cases}
\end{equation}
for $\ell=1,\dotsc,m$,
with $\theta_j^0$ satisfying \cref{PT2B}.
It also follows from this construction that
\cref{PT2A} has a solution for $(\cG,\Theta)$ if and only
if it is solvable for $(\cG^\ell,\Theta^\ell)$, and that
\begin{equation*}
\kappa_{ij}(\cG,\Theta) \,=\, \kappa_{ij}(\cG^\ell,\Theta^\ell)
\quad\forall\,(i,j)\in\cE^\ell\,,\quad \ell=0,\dotsc,m\,.
\end{equation*}

Therefore, since $\cI^{m}$ is a singleton,
say $\cI^{m}=\{\Hat\imath\}$,
\cref{PT2A} has a solution for $(\cG,\Theta)$ if and only if
\begin{equation}\label{PT2F}
\begin{aligned}
\vartheta_p\, p_{\Hat\imath}
&\,=\, -\Hat\lambda_{\Hat\imath} + \sum_{j\in\cJ(\Hat\imath)}
\mu_{\Hat\imath j}  \theta_j^{m}\\
&\,=\, -\Hat\lambda_{\Hat\imath} +
\sum_{j\in\cJ(\Hat\imath)} \dd(\Hat\imath,j) \theta_j^{m-1}
-\dd(\Hat\imath,\imath_{m})\, \Hat\lambda_{\imath_{m}}
 - \vartheta_p\, \dd(\Hat\imath,\imath_m)\,p_{\imath_m}
+ \sum_{k\in\cJ(\imath_{m})\setminus\{\jmath_{m}\}}
\dd(\Hat\imath,k)\,\theta_k^{m-1}\\
&\,=\, \sum_{j\in\cJ^{m-1}} \dd(\Hat\imath, j) \theta_j^{m-1}
-\sum_{i\in\cI^{m-1}} \dd(\Hat\imath,i) \Hat\lambda_i
- \vartheta_p \,\dd(\Hat\imath,\imath_m)p_{\imath_m}\,,
\end{aligned}
\end{equation}
where in the second equality we use \cref{PT2E}, and in the third
equality we use the fact that $\jmath_{m}\in\cJ^{m-1}$ which
is true by construction.
Next, an easy calculation using \cref{PT2E} shows that
\begin{equation}\label{PT2G}
\begin{aligned}
\sum_{j\in\cJ^{\ell}} \dd(\Hat\imath,j) \theta_j^{\ell}
-\sum_{i\in\cI^{\ell}} \dd(\Hat\imath,i)\Hat\lambda_i
&\,=\,
\sum_{j\in\cJ^{\ell-1}} \dd(\Hat\imath,j) \theta_j^{\ell-1}
-\sum_{i\in\cI^{\ell-1}} \dd(\Hat\imath,i)\Hat\lambda_i
- \vartheta_p\, \dd(\Hat\imath,\imath_\ell)p_{\imath_\ell}
\end{aligned}
\end{equation}
for $\ell =1,\dotsc,m$.
Therefore, using the recursion \cref{PT2G} in \cref{PT2F} we obtain
\begin{equation}\label{PT2H}
\begin{aligned}
\vartheta_p\, p_{\Hat\imath}
&\,=\, \sum_{j\in\cJ^0} \dd(\Hat\imath,j) \theta_j^{0}
-\sum_{i\in\cI^0} \dd(\Hat\imath,i)\Hat\lambda_i
- \vartheta_p \sum_{i\in\cI^0\setminus \{\Hat\imath\}}
\dd (\Hat\imath,i)p_i\\
&\,=\, \sum_{j\in\cJ} \dd(\Hat\imath,j) \theta_j
-\sum_{i\in\cI} \dd(\Hat\imath,i)\Hat\lambda_i
- \vartheta_p \sum_{i\in\cI\setminus \{\Hat\imath\}}
\dd(\Hat\imath,i) p_i\,,
\end{aligned}
\end{equation}
where in the last equality we use \cref{PT2B}.
Solving \cref{PT2H}, we obtain
\begin{equation}\label{PT2I}
\vartheta_p \,=\, \frac{1}{\sum_{\ell\in\cI}\dd(\Hat\imath,\ell)p_\ell}
\Biggl(\sum_{j\in\cJ} \dd(\Hat\imath,j)\,\theta_j
-\sum_{i\in\cI} \dd(\Hat\imath,i)\Hat\lambda_i\Biggr)\,.
\end{equation}
Note that the fractions
$\frac{\sum_{j\in\cJ} \dd(i,j)}{\sum_{\ell\in\cI}\dd(i,\ell)p_\ell}$ and
$\frac{\sum_{\ell\in\cI} \dd(i,\ell)}{\sum_{\ell\in\cI}\dd(i,\ell)p_\ell}$
do not depend on $i\in\cI$.
This can be seen, for example, by multiplying the numerator and denominator
by $\dd(i',i)$ and using the multiplicative property of the function $\dd$.
This fact together with \cref{PT2I} establishes \cref{ET2A}.
\Halmos\endproof

{\color{dred}
\begin{example}
To better illustrate the proof and the notations used, we show how the steps
in the proof are applied to the network in \cref{Fig1}. 
\begin{figure}[H]
\centerline{\includegraphics[height=1.8in]{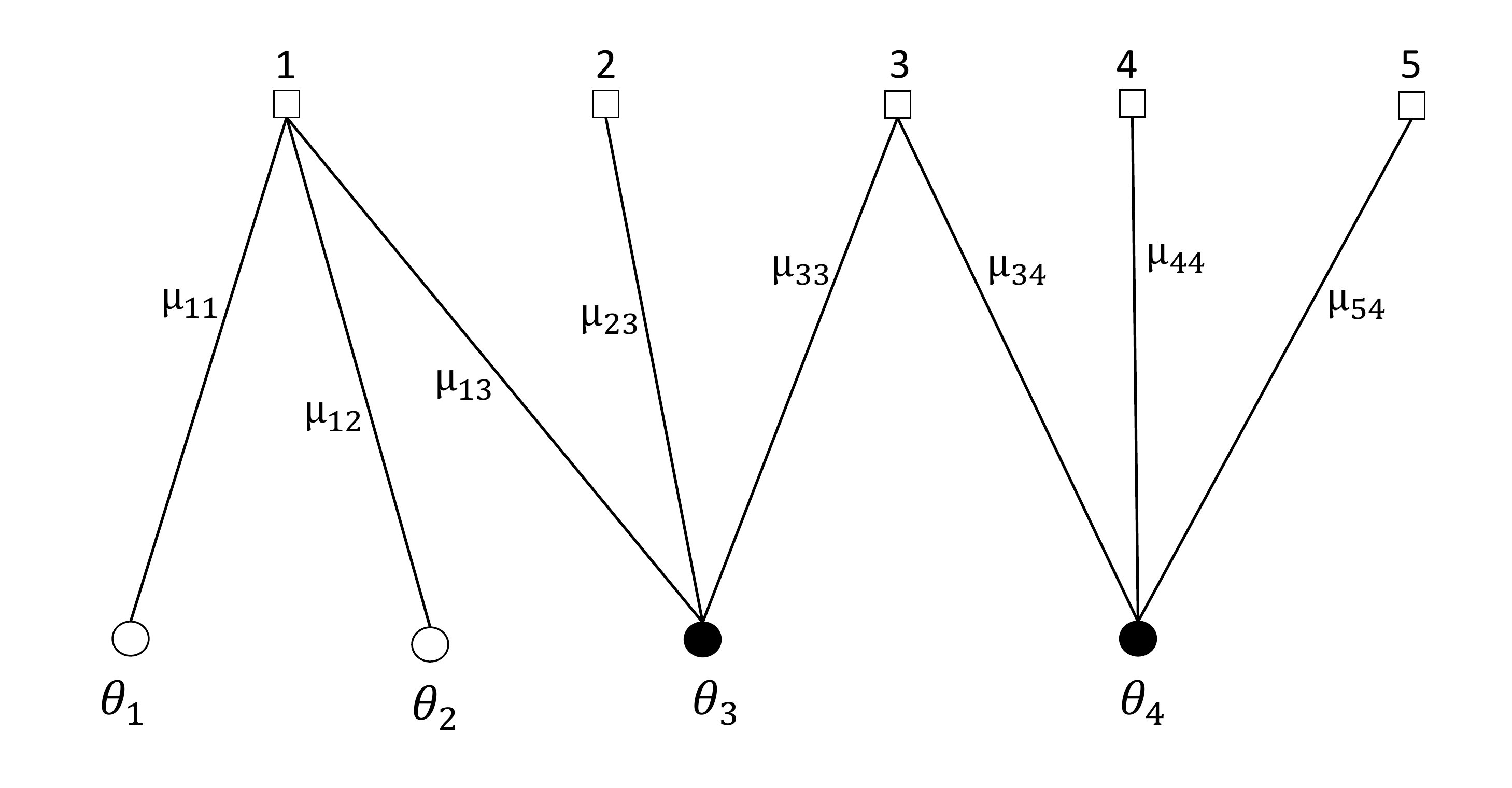}}
\centerline{\caption{A network with 5 classes  and 4 pools}}
\label{Fig1}
\end{figure}
Since $\cI_{\mathsf{leaf}}$ is the set of leaf classes, then $\cI_{\mathsf{leaf}} = \{2,4,5\}$
and $\widetilde{\cJ} = \{3,4\}$. These are the customer leaves to be removed in Step 1 of the algorithm and are shown as dashed edges.
\begin{figure}[H]
\centerline{\includegraphics[width=0.9\textwidth]{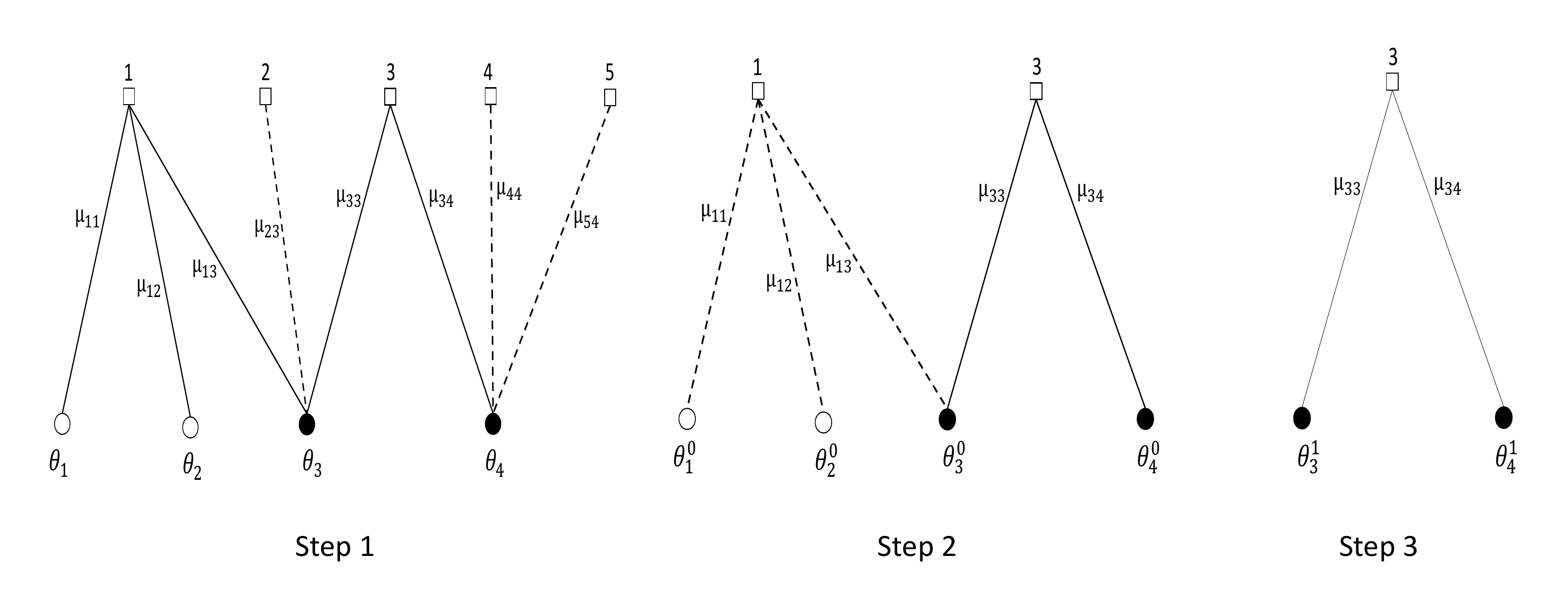}}
\centerline{\caption{Illustration of the algorithm}}
\end{figure}
In addition, the updates of the network parameters computed using \cref{PT2B} are given by
\begin{equation*}
\begin{pmatrix}
\theta_1^0\\[4pt]
\theta_2^0\\[4pt]
\theta_3^0\\[4pt]
\theta_4^0
\end{pmatrix} \,=\,
\begin{pmatrix}
\theta_1\\[4pt]
\theta_2\\[4pt]
\theta_3 - \mu_{23}^{-1}(\vartheta_p p_2 + \Hat{\lambda}_2)\\[4pt]
\theta^4 - \mu_{44}^{-1}(\vartheta_p p_4 + \Hat{\lambda}_4)- \mu_{54}^{-1}(\vartheta_p p_5 + \Hat{\lambda}_5)
\end{pmatrix}\,.
\end{equation*}
In Step 2, class $1$ is selected i.e., $\imath_1 = 1$ and $\jmath_1 = 3$. The parameters are updated according to \cref{PT2C} and are given by
\begin{equation*}
\begin{pmatrix}
\theta_3^1\\[4pt]
\theta_4^1
\end{pmatrix} \,=\, 
\begin{pmatrix}
\theta_3^0 - \mu_{13}^{-1}\bigl(\Hat{\lambda}_1 + \vartheta_p p_1 - \mu_{11}\theta_1^0 - \mu_{12}\theta_2^0\bigr)\\[4pt]
\theta_4^0
\end{pmatrix}\,,
\end{equation*}
and class $1$ is removed
together with the associated edges.
The resulting network is shown in Step~3 and has only class $3$, i.e., $\Hat{\imath} = 3$.
Therefore, the SWSS $\vartheta_p$ is computed using \cref{PT2F} as
\begin{equation*}
\vartheta_p \,=\, -\Hat{\lambda}_3 + \mu_{33} \theta_3^1 + \mu_{34}\theta_4^1\,.
\end{equation*}
\end{example}
}

{\color{dred}
\begin{remark}\label{R3}
In analogy to the definitions in the beginning of the section,
we define the gain
$\dd(j,j')$
between any pair $j,j'\in \cJ$, $j\ne j'$,
by
\begin{equation*}
\dd(j,j') \,\df\,
\prod_{k=1}^{m-1} \frac{\mu_{i_{k+1}j_{k+1}}}{\mu_{i_{k+1}j_k}}\,,
\end{equation*}
where the product is evaluated over the
analogous path  $(j_1, i_2, j_2,\dotsc, i_m, j_m)$
connecting $j\equiv j_1$ to $j'\equiv j_m$ in $\cG$.
Suppose for simplicity that $\lambda_i^n=n\lambda_i$ and $\mu_{ij}^n=\mu_{ij}$.
It follows from \cref{ET2A} that if we decrease $\Hat\nu_j$ by an amount $\delta$
and increase $\Hat\nu_{j'}$ by an amount $\dd(j,j')\delta$, then the value of
$\vartheta_p$ stays the same.
This has the following interpretation.
The contribution of one server at pool $j$ in the stability
of the network is the same as that of $\dd(j,j')$ servers
in pool $j'$.
\end{remark}
}

With $p$ a positive vector in $\varDelta_I$,  we define
\begin{equation*}
R_i \,\df\, \frac{\sum_{j\in\cJ} \dd(i,j) \theta_j
- \sum_{\ell\in\cI} \dd(i,\ell) \Hat\lambda_\ell}
{p_i}\,,\quad\text{and} \quad
\varGamma_i \,\df\,
\sum_{\ell\in\cI}\dd(i,\ell) \frac{p_\ell}{p_i}
\end{equation*}
for $i\in\cI$. Note that 
\begin{equation}\label{EC1B}
\sum_{i\in\cI} \varGamma_i^{-1}
\,=\, \sum_{i\in\cI} \frac{p_i}{\sum_{\ell\in\cI}\dd(i,\ell)p_\ell}
\,=\, \sum_{i\in\cI}\frac{\dd(1,i)p_i}{\sum_{\ell\in\cI}\dd(1,\ell)p_\ell}\,=\,1\,.
\end{equation}
By \cref{T2} we have $\varGamma_i^{-1}\vartheta_p^{-1}
= \frac{1}{R_i}$ for all $i\in\cI$,
and summing up this equality over $i\in\cI$, and using \cref{EC1B},
we obtain the following corollary.

\begin{corollary}\label{C1}
It holds that
\begin{equation*}
\frac{1}{\vartheta_p}\,=\,\sum_{i\in\cI} \frac{1}{R_i}\,.
\end{equation*}
\end{corollary}

As asserted in \cref{C1} the SWSS parameter $\vartheta_p$ is the harmonic mean of
the variables $R_i$.
This could be compared with the formula of the resistance of branches connected
in parallel in electric circuits.

\begin{remark}
With $e_i\in\RR^I$  as defined in \cref{S1.1}, we have the identity
$\vartheta_{e_i}=R_i$ for all $i\in\cI$.
In other words, $R_i$ is the maximum permissible safety staffing for class $i$
without allowing the safety staffing of the other classes to go negative.
\end{remark}

\section{Relating the SWSS to the Drift of the Diffusion Limit}\label{S4}

\rd{In this section, we establish a characterization of the SWSS
$\vartheta_p$ in terms of the parameters of the diffusion limit of the Markovian network.
We also obtain an analogous characterization of $\vartheta_p^n$.
The key results are in \cref{P1,P2,T3,T4}, and these are essential for the
stability analysis in \cref{S5,S6}.}

For each $i\in\cI$ and $j\in\cJ$, we let $X_i^n=\{X_i^n(t)\colon t\ge 0\}$
denote the total number of class $i$ customers in the system
(both in service and in queue), $Z_{ij}^n=\{Z_{ij}^n(t),\,t\ge0 \}$ the number of
class $i$ customers currently being served in pool $j$,
$Q_i^n=\{Q_i^n(t),\,t\ge0\}$ the number of class $i$ customers in the queue,
and $Y_j^n=\{Y_j^n(t),\,t\ge0\}$ the number of idle servers in server pool $j$.
Let $X^{n} = (X_i^{n})_{i \in \cI}$, $Y^{n} = (Y_j^{n})_{j \in \cJ}$, 
$Q^{n} = (Q_i^{n})_{i \in \cI}$,
and $Z^{n} = (Z_{ij}^{n})_{(i,j)\in\cE}$.
The process $Z^{n}$ is the scheduling control.
Let $(x,z)\in\ZZ_+^I\times \ZZ^{\cG}_+$ denote a state-action pair.
We define
\begin{equation*}
q_i(x,z) \,\df\, x_i - \sum_{j\in\cJ} z_{ij}\,,\quad i\in\cI\,,\qquad
y_j^n(z) \,\df\,   N_j^n - \sum_{i\in\cJ}z_{ij}\,,
\quad j\in\cJ\,,
\end{equation*}
and the (work-conserving) \emph{action space} $\cZn(x)$ by
\begin{equation*}
\cZn(x)\,\df\, \bigl\{z \in \ZZ^{\cG}_+\,\colon
q_i(x,z) \wedge y_j^n(z) =0\,,~
q_i(x,z)\ge0\,,~y_j^n(z) \ge0 \quad\forall\,(i,j)\in\cE\bigr\}\,. 
\end{equation*}

\subsection{The Diffusion Scaling}\label{S4.1}
\rd{In this section, we first write the infinitesimal generator of the state process
for the $n^{\mathrm th}$ system in the form of \cref{E-brLgn}
which is a linear expression involving second order
and first order difference quotients.
The coefficients of the first order difference quotients comprise the
`drift' of the infinitesimal generator.
Then, in \cref{P1}, we obtain an explicit form of the drift.
The components of the constant term in the drift, which are denoted as $h^n_i$,
are given in \cref{E-hn}, and are equal
to the diffusion-scaled difference between the service rates $\lambda_i^n$
and the service rate of class $i$ when the service allocations
are chosen as the solution of \cref{E-LP}.
The constant term $h^n$ and the matrices $B_1^n$ and $B_2^n$ that describe the
linear part of the drift, are used to characterize $\vartheta_p^n$ in
\cref{T4} in \cref{S4.2}.}

We introduce some suitable notation to describe the diffusion scale.
For additional details see \citet{AP19}.
With $\xi^*\in\RR^\cG_+$ being the solution of the \cref{E-LP}, we
define  $\Bar{z}^n \in\RR^{\cG}_+$ and $\Bar{x}^n\in\RR^I$ by
\begin{equation}\label{E4.1A}
 \Bar{z}^n_{ij} \,\df\, \xi^*_{ij}N_j^n\,,
 \qquad \Bar{x}^n_i \,\df\, \sum_{j\in\cJ} \Bar{z}^n_{ij}\,,
\end{equation}
and
\begin{equation}\label{E4.1B}
\Breve{x}^{n}\,=\, \Breve{x}^{n}(x) \,\df\, \frac{x -  \Bar{x}^n}{\sqrt{n}}\,, \qquad
 \Breve{z}^n \,=\, \Breve{z}^{n}(z)\,\df\, \frac{z -  \Bar{z}^n}{\sqrt{n}}\,,
\end{equation}
for $x\in\ZZ^{I}_{+}$ and $z\in \cZn(x)$.
We also let $\sS^n$ denote the state space in the diffusion scale, that is,
$\sS^n \df
\bigl\{\Breve{x}\in\Rm\,\colon \sqrt{n}\Breve{x}+\Bar{x}^n\in\ZZ^I_{+}\bigr\}$.
The diffusion-scaled variables are defined by
\begin{equation} \label{E4.1C}
\Breve{X}^{n}_i(t) \,\df\,\Breve{x}^{n}\bigl(X_i^{n}(t)\bigr)\,,  \quad
\Breve{Z}^{n}_{ij}(t) \,\df\,\Breve{z}^{n}\bigl(Z_{ij}^{n}(t)\bigr)\,,\quad
\Breve{Q}^{n}_i(t) \,\df\, \frac{Q_i^{n}(t)}{\sqrt{n}}  \,,
\quad\text{and} \quad
\Breve{Y}^{n}_j(t) \,\df\, \frac{Y_j^{n}(t)}{\sqrt{n}}  \,.
\end{equation}
Under a stationary Markov policy
$Z^n(t) = z(X^n(t))$ for some function $z\colon\ZZ^{I}_+ \to \ZZ^{\cG}_+$,
the process $X^n$ is Markov with controlled generator 
\begin{equation}\label{E-Lgn}
\Lg_z^n f(x) \,\df\, \sum_{i\in\cI}\Biggl(\lambda_i^n\bigl(f(x+e_i)-f(x)\bigr)
+ \sum_{j\in\cJ(i)}\mu_{ij}^n z_{ij}\bigl(f(x-e_i)-f(x)\bigr)\Biggr)
\end{equation} 
for $f\in C(\RR^I)$ and $x \in \ZZ^{I}_{+}$.

We drop the explicit dependence on $n$ in the diffusion-scaled
variables in order to simplify the notation.
Note that a work-conserving stationary Markov policy $z$,
that is a map $z\colon\ZZ^{I}_+ \to \ZZ^{\cG}_+$
such that $z(x)\in\cZn(x)$ for all $x\in\ZZ_+^I$, gives rise to a
stationary Markov policy
$\Breve{z}\colon\sS^n\to\RR^\cG$, with
\begin{equation*}
\Breve{z}(\Breve{x}) \,\in\,\bcZn(\Breve{x}) \,\df\,
\bigl\{\Breve{z}\,\colon
\sqrt{n}\Breve{z}+ \Bar{z}^n\in\cZn(\sqrt{n}\Breve{x}+ \Bar{x}^n)\bigr\}
\qquad \forall\,\Breve{x}\in\sS^n\,,
\end{equation*}
via \cref{E4.1B} (and vice-versa).
Let $h^{n} = (h^{n}_1,\dotsc,h^{n}_I)\transp$ be defined by
\begin{equation}\label{E-hn}
h^{n}_i\,\df\, \frac{1}{\sqrt{n}} \Biggl(\lambda^{n}_i
-  \sum_{j \in \cJ(i)} \mu_{ij}^{n} \xi^*_{ij}N_j^n \Biggr) \,.
\end{equation}
By  the assumptions on the parameters in
\cref{E-Par01,E-Par02}, we have $h^{n}_i\to h_i$ as $n\to\infty$,
with $h_i$ as defined in \cref{E-h}.
We let $h\df(h_1,\dotsc,h_I)\transp$.
Using \cref{E4.1A,E4.1B,E-Lgn,E-hn} and rearranging terms, the controlled generator
of the corresponding
diffusion-scaled process can be written as
\begin{equation}\label{E-brLgn}
\begin{aligned}
\Breve\Lg_{\Breve{z}}^n f(\Breve{x})\,\df\, &\,=\,
\sum_{i\in\cI}\frac{\lambda_i^n}{n}\,
\frac{ f\bigl(\Breve{x}+\tfrac{1}{\sqrt n} e_i\bigr) - 2f(\Breve x)
+ f\bigl(\Breve{x}-\tfrac{1}{\sqrt n} e_i\bigr)}{n^{-1}} \\
&\mspace{50mu}-\sum_{i\in\cI} b^n_i(\Breve{x},\Breve{z})\,
\frac{f\bigl(\Breve{x} -\tfrac{1}{\sqrt n} e_i\bigr)
- f(\Breve x)}{n^{-\nicefrac{1}{2}}}\,,
\qquad \Breve{x}\in\sS^n\,,\ \Breve{z}\in\bcZn(\Breve{x})\,, 
\end{aligned}
\end{equation}
where the `drift' $b^n=(b_1^n,\cdots,b_I^n)\transp$ is given by
\begin{equation}\label{E-GnC}
b^n_i(\Breve{x},\Breve{z}) \,\df\, h_i^n - \sum_{j\in\cJ(i)}\mu_{ij}^n
\Breve{z}_{ij}\,,\quad\Breve{z}\in \bcZn(\Breve{x})\,,\ i\in\cI\,.
\end{equation}

\begin{definition}\label{D2}
For $\Breve{x}\in\sS^n$ and $\Breve{z}\in\bcZn(\Breve{x})$, we define
\begin{equation}\label{ED2A}
\Breve{q}^n_i(\Breve{x},\Breve{z})\,\df\, \Breve{x}_i
- \sum_{j\in\cJ(i)} \Breve{z}_{ij}\,,\quad i\in\cI\,,
\qquad
\Breve{y}^n_j(\Breve{z})\,\df\, - \sum_{i\in\cI(j)} \Breve{z}_{ij}\,, \quad j\in\cJ\,,
\end{equation}
and
$\Breve{\zeta}^n(\Breve{x},\Breve{z})\df
\langle e, \Breve{q}^n(\Breve{x},\Breve{z})\bigr\rangle\wedge\,
\langle e, \Breve{y}^n(\Breve{z})\bigr\rangle$.
We also let
\begin{equation}\label{E-Act}
\Act \,\df\, \varDelta_I \times \varDelta_J \,\df\, \bigl\{ u^c \in \RR^{I}_+
\,\colon\, \langle e, u^c \rangle \,=\, 1\bigr\} \times \bigl\{u^s \in \RR^J_+
\,\colon\,\langle e, u^s \rangle \,=\, 1\bigr\}\,.
\end{equation}
\end{definition}

In the following proposition, we give an explicit expression of the drift $b^n$ and the some of its structural properties. 

{\color{dred}
\begin{proposition}\label{P1}
For any $\Breve{z}\in\bcZn(\breve{x})$ with $\Breve{x}\in\sS^n$,
there exists $u=u(\Breve{x},\Breve{z})\in\Act$ such that the drift
$b^n$ in \cref{E-GnC} takes the form
\begin{equation}\label{E-bn}
b^n(\Breve{x},\Breve{z}) \,=\, h^n 
- B_1^n \bigl(\Breve{x}-\langle e,\Breve{x}\rangle^{+} u^c\bigr)
+ B_2^n u^s \langle e,\Breve{x}\rangle^{-}
+\Breve{\zeta}^{n}(\Breve{x},\Breve{z}) \bigl(B_1^n u^c + B_2^n u^s\bigr)\,.
\end{equation}
In \cref{E-bn}, $h^n = (h^n_1, \dotsc, h^n_I)\transp$ with $h^n_i$ is as in \cref{E-hn},
$B_1^n\in\RR^{I\times I}$, and $B_2^n\in\RR^{I\times J}$.

In addition, given any $(\Hat\imath,\Hat\jmath)\in\cE$, there exists an ordered basis
$\mathfrak{D} =\bigl(\alpha, (\beta)_{-\Hat\jmath}\bigr)$,
with $\alpha_{\Hat\imath}$ being the last element of $\alpha$,
with respect to which the matrices $B_1^n$ and $B_2^n$
take the following form:
\begin{itemize}
\item[\ttup a]
    $B^n_1$ is  a lower-diagonal $I\times I$ matrix with positive diagonal
elements and $(B_1)_{II} = \mu^n_{\Hat\imath\Hat\jmath}$;
\item[\ttup b]
$B^n_2$ is an $I \times J$ matrix whose last column is identically zero.
\end{itemize}
\end{proposition}
}

\proof{Proof.}
Using \cref{D2}, it is easy to see that there exists $u=(u^c,u^s)\in\Act$,
depending on $\Breve{x}\in\sS^n$ and $\Breve{z}\in\bcZn(\Breve{x})$,
such that
\begin{equation}\label{ES4.1B}
\Breve{q}^n(\Breve{x},\Breve{z}) \,=\,\bigl(\Breve{\zeta}^n(\Breve{x},\Breve{z})
+\langle e, \Breve{x}\rangle^+\bigr)\,u^c\,,
\quad\text{and}\quad
\Breve{y}^{n}(\Breve{z}) \,=\, \bigl(\Breve{\zeta}^n(\Breve{x},\Breve{z})
+\langle e, \Breve{x}\rangle^-\bigr)\,u^s\,.
\end{equation}
Let $\sD\df\bigl\{ (\alpha, \beta) \in \RR^{I}\times \RR^{J}\,\colon\,
\langle e, \alpha \rangle \,=\, \langle e, \beta \rangle\bigr\}$.
Define the
linear map $\Psi=[\Psi_{ij}]\colon \sD\to \RR^{I\times J}$ as the solution of 
\begin{equation} \label{E-Psi}
\sum_j \Psi_{ij}(\alpha,\beta)\,=\,\alpha_i \quad \forall i \in \cI\,,
\quad\text{and}\quad
\sum_{i} \Psi_{ij}(\alpha,\beta)\,=\,\beta_j \quad \forall j \in \cJ\,,
\end{equation}
with $\Psi_{ij}(\alpha,\beta) = 0$ for $i \nsim j$.
It is shown in Proposition~A.2 of \citet{Atar-05a} that 
if $\cG$ is a tree, the linear map $\Psi$ is unique. 
Since
$\bigl(\Breve{x}-\Breve{q}^{n}(\Breve{x},\Breve{z}),- \Breve{y}^{n}(\Breve{z})\bigr)
\in \sD$
by \cref{ED2A},
using the linearity of the map $\Psi$ and \cref{ES4.1B,E-Psi}, it follows that
\begin{equation*}
\begin{aligned}
\Breve{z} &\,=\,
\Psi\bigl(\Breve{x}-\Breve{q}^{n}(\Breve{x},\Breve{z}),- \Breve{y}^{n}(\Breve{z})\bigr)\\
&\,=\, \Psi\bigl(\Breve{x}-\langle e,\Breve{x}\rangle^{+} u^c,
-\langle e,\Breve{x}\rangle^{-} u^s\bigr)
-\Breve{\zeta}^{n}(\Breve{x},\Breve{z})\,\Psi(u^c,u^s)\,.
\end{aligned}
\end{equation*}

Consider the matrices $B_1^n\in\RR^{I\times I}$ and
$B_2^n\in\RR^{I\times J}$ defined by
\begin{equation}\label{E-linear}
\sum_{j\in\cJ(i)} \mu^n_{ij} \Psi_{ij}(\alpha,\beta) \,=\,
\bigl(B_1^n \alpha + B_2^n \beta\bigr)_i\,,\quad
\forall\,i\in\cI\,,\ \forall (\alpha,\beta)\in \sD\,.
\end{equation}
As shown in \citet[Lemma~4.3]{AP16},  given any $(\Hat\imath,\Hat\jmath)\in\cE$,
there exists an ordered basis
$\mathfrak{D} =\bigl(\alpha, (\beta)_{-\Hat\jmath}\bigr)$, $j\in\cJ$, where
$(\beta)_{-\Hat\jmath}= \{\beta_\ell\,,\, \ell\ne \Hat\jmath\}$
and $\alpha_\imath$ is the last element of $\alpha$, such
that the matrix $B^n_1$ in \cref{E-linear} satisfies assertion (a).
It is also clear from the proof of the above lemma, that $B^n_2$
satisfies assertion (b).
This completes the proof.
\Halmos\endproof


\subsection{The Diffusion Limit}\label{S4.2}

\rd{In this section we present some important properties of the drift of
the limiting diffusion.
This takes the form of the piecewise-affine function given in \cref{E-drift2}.
The coefficients $h$, $B_1$, and $B_2$ are the limits of
$h^n$, $B_1^n$, and $B_2^n$ in \cref{E-bn}, respectively, as $n\to\infty$.
Two additional important results are presented in this section:
\Cref{T3} which characterizes the gains $\dd(i,\ell)$ for $i,\ell \in \cI$ in terms of the matrix
$B_1$, and \cref{T4} which obtains an analogous characterization of
the SWSS parameters $\vartheta_p$ and $\vartheta^n_p$.} 


To discuss the diffusion limit we need the concept of \emph{joint work conservation}.
We say that an action $\Breve{z}\in\bcZn(\Breve{x})$ is \emph{jointly work
conserving} (JWC), if $\Breve{\zeta}^n(\Breve{x},\Breve{z})=0$ \rd{, i.e., a scheduling rule
for which either there are no customers in the queues, or no server in the system is idle.
Simple examples show that it is not possible to specify such an action on the
whole state space.}
However, as shown in \citet[Lemma~3]{Atar-05b}, there exists
$M_0>0$ such that the collection
of sets $\Breve{\sX}^{n}$ defined by
\begin{equation}\label{E-jwc}
\Breve{\sX}^{n}\,\df\,\bigl\{\Breve{x} \in
\sS^n\,\colon  \norm{\Breve{x}}^{}_1\le M_{0}\, \sqrt{n}  \bigr\}
\end{equation}
has the following property:
For any $\Breve{x}\in\Breve{\sX}^{n}$ and a pair $(\Breve{q},\Breve{y})$
such that $\sqrt{n}\Breve{q}\in\ZZ^I_+$, $\sqrt{n}\Breve{y}\in\ZZ^J_+$, and
\begin{equation*}
\langle e,\Breve{q}\rangle\wedge \langle e,\Breve{y}\rangle \,=\,0\,,\quad
\langle e,\Breve{x}-\Breve{q}\rangle \,=\, \langle e,-\Breve{y}\rangle\,,
\quad\text{and}\quad
\Breve{y}_j\sqrt{n}\,\le\, N^n_j\,,\quad j\in\cJ\,,
\end{equation*}
 it holds that
$\Psi(\Breve{x}-\Breve{q},-\Breve{y})\in\bcZn(\Breve{x})$.

Under any stationary Markov scheduling policy that
is jointly work-conserving
in the set $\Breve{\sX}^{n}$, the diffusion-scaled state process $\Breve{X}^n$ converges
weakly to a limit $X$ described as follows. 
For $u\in\Act$, let $\Breve{\Psi}[u]\colon\RR^I\to\RR^{\cG}$ be defined by
\begin{equation}\label{E-bPsi}
\Breve{\Psi}[u](x) \,\df\,
\Psi(x- \langle e, x\rangle^{+} u^c, - \langle e, x\rangle^{-} u^s)\,,
\end{equation}
where $\Psi$ is as in \cref{E-Psi}.
The limiting controlled diffusion $ X$ is given by the It\^o equation 
\begin{equation} \label{E-diff}
d X_t \,=\, b( X_t, U_t)\,\D{t} + \Sigma \, \D W_t\,, 
\end{equation}
where $W$ is an $I$-dimensional standard Wiener process, and 
$\Sigma \df \diag\bigl(\sqrt{2 \lambda_1}, \dotsc, \sqrt{2 \lambda_I}\bigr)$.
The drift $b\colon \RR^I \times \Act \to \RR^I$ takes the form
\begin{equation}\label{E-drift1}
 b_i(x,u)\,=\, b_i\bigl( x, (u^c, u^s)\bigr) \,\df\, h_i - \sum_{j \in \cJ(i)}
\mu_{ij} \Breve{\Psi}_{ij}[u](x) \qquad\forall\,i\in\cI\,,
\end{equation}
where $\Breve{\Psi}_{ij}[u]$ is as in \cref{E-bPsi} and
$h_i$ is given by \cref{E-h}.
This result was first shown in \citet{Atar-05a,Atar-05b}.
We focus on the class $\Usm$ of stationary Markov
controls, that is, $U_t=v(X_t)$ for some measurable function $v\colon \RR^I\to\Act$.

The following proposition is the exact analog of \cref{P1},
and follows by taking limits as $n\to\infty$ in \cref{E-bn}
and employing the convergence of the parameters in \cref{E-Par01,E-Par02}.

\begin{proposition}\label{P2}
Given any $(\Hat\imath,\Hat\jmath)\in\cE$, there exists an ordered basis
$\mathfrak{D} =\bigl(\alpha, (\beta)_{-\Hat\jmath}\bigr)$,
with $\alpha_{\Hat\imath}$ being the last element of $\alpha$, with respect to which
the drift $b$ in \cref{E-drift1} has the representation
\begin{equation} \label{E-drift2}
b(x, u) \,=\, h - B_1 (x -  \langle e, x \rangle^{+} u^c)
+ \langle e, x \rangle^{-} B_2 u^s\,,
\end{equation}
where
\begin{itemize}
\item[\ttup a]
$h = (h_1, \dotsc, h_I)\transp$ and $h_i$ is as in \cref{E-h},
\item[\ttup b]
$B_1$ is  a lower-diagonal $I\times I$ matrix with positive diagonal
elements and $(B_1)_{II} = \mu_{\Hat\imath\Hat\jmath}$\,, 
\item[\ttup c]
$B_2$ is an $I \times J$ matrix whose last column is identically zero.
\end{itemize}
\end{proposition}

For $f\in C^2(\RR^m)$, we define
\begin{equation}\label{E-Lg}
\Lg_u f(x) \,\df\, \frac{1}{2} \trace\bigl(\Sigma\Sigma\transp \nabla^2f(x)\bigr)
+\bigl\langle b(x,u),\nabla f(x)\bigr\rangle\,,
\end{equation}
with $\nabla^2f$ denoting the Hessian of $f$.

\begin{remark}\label{R5}
We remark that \cref{E4.1C} differs from the 
usual definition of the
diffusion-scaled processes found in the literature
(see \citet{Atar-05b} and \citet{AP16,AP18,AP19}).
We refer to the processes $\Breve{X}^n_i$ as the ``re-centered''
diffusion-scaled processes.
One may also center the process $X^n_i$ around $n x^*_i$  where 
$x^*$ is defined in  \cref{E-statfl}.
It is clear that the limit processes using these different centering terms
only differ in the drift by a constant,
and therefore they are equivalent as far as their ergodic properties
are concerned. See also \cref{R8}.
\end{remark}

\begin{remark}\label{R6}
Note that if the networks have renewal arrivals and the
service times are exponential, we again obtain a diffusion 
limit for the above diffusion-scaled processes, which has 
the same drift as the Markovian case, and whose covariance matrix captures
the variability in the arrivals processes. In 
particular, if the class-$i$ arrival process $A^n_i$ is 
renewal with interarrival times of rate $\lambda^n_i$ 
(satisfying \cref{E-Par01,E-Par02}) and 
variance $(\sigma^n_i)^2$ (satisfying $\sigma^n_i \to \sigma_i>0$ as $n\to \infty$), then
\begin{equation*}
\Hat{A}^n_i(t) \,=\,
n^{\nicefrac{1}{2}} (A^n_i(t) - \lambda^n_i t) \,\Rightarrow\,
\Hat{A}_i(t) \,=\, W_i(\lambda_i c^2_{a,i}t)\,,
\end{equation*}
where $W_i$ is a
standard Brownian motion, and
$c^2_{a,i} \df \lambda_i^2 \sigma_i^2$.
As a consequence, the covariance matrix $\Sigma$ in \cref{E-diff} takes
the form 
$
\Sigma= \diag\bigl(\bigl( \lambda_1(1+c_{a,1}^2)\bigr)^{\nicefrac{1}{2}}, \dotsc, \bigl( \lambda_I(1+c_{a,I}^2) \bigr)^{\nicefrac{1}{2}}\bigr)\,.
$
Thus, the results in \cref{T3} also hold for the 
networks with renewal arrivals and exponential service times.
The same applies to the results regarding the limiting controlled diffusion
in \cref{T5,T7}. See also \cref{R9} for the results concerning the 
diffusion-scaled processes.
\end{remark}

The following result is essential in proving the main theorem of this section.
Recall the definition in \cref{E-dd}.

\begin{theorem}\label{T3}
It holds that
\begin{equation}\label{ET3A}
\dd(i,\ell) \,=\, \frac{\bigl(e\transp B_1^{-1}\bigr)_{\ell}}
{\bigl(e\transp B_1^{-1}\bigr)_{i}}\qquad\forall\, i,\ell\in\cI\,.
\end{equation}
In addition, $\bigl(e\transp B_1^{-1}\bigr)_i>0$ for all $i\in\cI$.
\end{theorem}
\proof{Proof.}
We start with the following observation: 
Let $\vartheta_p$,  $\kappa=[\kappa_{ij}]\in\RR^{\cG}$ be the unique solution of \cref{E-LP'}. This means that 
$\sum_{j \in \cJ(i)}\mu_{ij} \kappa_{ij} = \hat{\lambda}_i + \vartheta_p p_i$
and
$\sum_{i \in \cI(j)} \kappa_{ij} = \theta_j$ for $i \in \cI$ and $j \in \cJ$. Let $\alpha_i = \sum_{j} \kappa_{ij}$ and note that $e\transp \alpha = e\transp \theta$ where $\alpha = (\alpha_i)_{i \in \cI}$ and $\theta = (\theta_j)_{j \in \cJ}$. 

Using a similar approach to \cref{E-Psi,E-linear}, it follows that $\kappa_{ij} = \Psi_{ij}(\alpha, \theta)$ and $B_1\alpha + B_2\theta = \sum_{j \in \cJ(i)}\mu_{ij}\Psi_{ij}(\alpha, \theta) = \hat{\lambda} + \vartheta_p p$. Inverting $B_1$ and using that $e^T\alpha = e^T \theta$, one reaches
that 
\begin{equation}\label{EPT3A}
\vartheta_p e^T B_1^{-1} p  = e^T\theta + e^T B_1^{-1} (- \hat{\lambda} + B_2 \theta).
\end{equation}
This means that $\vartheta_p e^TB_1^{-1}p$ remains
constant when $p$ varies. Using \cref{ET2A,EPT3A} we get that $e^T B_1^{-1} p  = c \sum_{\ell \in \cI} \dd(i,\ell) p_\ell$ for every positive probability vector $p$ where $c$ is a constant independent of $p$. This of course implies that $\bigl(e^T B_1^{-1}\bigr)_\ell  = c\,\dd(i,\ell)$ for $\ell \in \cI$. It remains to note that $\dd(i,i) = 1$ by definition, see \cref{E-dd}, to conclude that $c = \bigl(e^T B_1^{-1}\bigr)_i$.

To prove the last assertion of the theorem, note that
$\bigl(e\transp B_1^{-1}\bigr)_I>0$ since $B_1$ is a lower diagonal matrix
with positive elements.
Thus, using \cref{ET3A}, we obtain
$\bigl(e\transp B_1^{-1}\bigr)_i = \dd(I,i) \bigl(e\transp B_1^{-1}\bigr)_I>0$ 
and
this concludes the proof.
\Halmos\endproof

\begin{theorem}\label{T4}
The variables $\vartheta_p$  and $\vartheta^n_p$ 
in \cref{D1} satisfy
\begin{equation} \label{ET4A}
\vartheta_p \,=\,  -\frac{\langle e, B_1^{-1} h \rangle}{\langle e, B_1^{-1} p\rangle}\quad
\text{and\ \ }
 \vartheta^n_p \,=\,  -\frac{\langle e, (B_1^n)^{-1} h^n\rangle}
{\langle e, (B_1^n)^{-1} p\rangle}\,,
\end{equation}
where $h=(h_i)_{i\in\cI}$ is given by \cref{E-h}, $h^n$ is defined in \cref{E-hn},
and $p\in\varDelta_I$ is a positive vector.
\end{theorem}

\proof{Proof.}
By \cref{T3}, we have
\begin{equation*}
e\transp B_1^{-1} p \,=\, \sum_{\ell\in\cI} \bigl(e\transp B_1^{-1}\bigr)_\ell\, p_\ell
\,=\, \bigl(e\transp B_1^{-1}\bigr)_i \sum_{\ell\in\cI} \dd(i,\ell)p_\ell\,,
\end{equation*}
and 
$e\transp B_1^{-1} h = \bigl(e\transp B_1^{-1}\bigr)_i \sum_{\ell\in\cI}
\dd(i,\ell)h_\ell$.
Combining these we obtain
\begin{equation*}
\begin{aligned}
-\frac{e\transp B_1^{-1} h}{e\transp B_1^{-1}p}
&\,=\,  \frac{\sum_{\ell\in\cI}\dd(i,\ell)\sum_{j\in\cJ(\ell)}\Bigl(\mu_{\ell j}
\xi^*_{\ell j}\Hat{\nu}_{\ell j}
+ \Hat{\mu}_{\ell j}z^*_{\ell j}\Bigr)
-\sum_{\ell\in\cI}\dd(i,\ell)\Hat{\lambda}_\ell}{\sum_{\ell\in\cI}\dd(i,\ell)p_\ell}
\,=\, \vartheta_p\,,
\end{aligned}
\end{equation*}
where in the last equality we used $\sum_{i\in\cI(j)}\xi^*_{ij} = 1$.
The same approach is used for $\vartheta^n_p$, thus establishing \cref{ET4A}.
\Halmos\endproof

\subsection{Some Examples}
In this part, we present some applications of \cref{T2,T4} by computing explicitly
the SWSS parameter for some networks.
We also give simple interpretations in the special case when
$\lambda_i^n = n\lambda_i$ and $\mu^n_{ij} = \mu_{ij}$, or equivalently, if
$\Hat{\lambda}_i = 0$, $\Hat\mu_{ij} = 0$, and $\theta_j = \Hat\nu_j$
for all $i\in\cI$ and $j\in\cJ$.

\begin{example}[The `N' Network]
For this network the SWSS parameter is given by
\begin{equation*}
\vartheta_p \,=\, \frac{\mu_{22}}{ \mu_{12}\,p_2 + \mu_{22} \, p_1}
\biggl(\mu_{11}\theta_1 + \mu_{12} \theta_2-\Hat\lambda_1 -\frac{\mu_{12}}{\mu_{22}}
\Hat\lambda_2\biggr)\,,
\end{equation*}
with 
\begin{equation*}
\kappa_{11}^* \,=\, \theta_1\,, \quad \kappa_{22}^*
\,=\, \frac{\Hat{\lambda}_2 + \vartheta_p\,p_2}{\mu_{22}}\,, 
\quad \kappa_{12}^* \,=\, \theta_2 - \kappa_{22}^*\,.   
\end{equation*}
In this case, the $B_1$ matrix is given by $B_1 = \diag (\mu_{12},\mu_{22})$
and the vector $h$ is given by 
\begin{equation*}
h \,=\, \begin{pmatrix}
\Hat{\lambda}_1 - \mu_{11}\xi^*_{11}\Hat{\nu}_1 - \Hat{\mu}_{11} z^*_{11}
- \mu_{12}\xi^*_{12}\Hat{\nu}_2 - \Hat{\mu}_{12} z^*_{12}\\[3pt]
\Hat{\lambda}_2 - \mu_{22}\xi^*_{22}\Hat{\nu}_2 - \Hat{\mu}_{22} z^*_{22}
\end{pmatrix}\,,
\end{equation*}
where $\xi^*_{11} = 1$, $\xi^*_{12} + \xi^*_{22} = 1$.
A simple calculation confirms that
$\vartheta_p = -\frac{\langle e,B_1^{-1} h \rangle}{\langle e,B_1^{-1}p\rangle}$.
In the special case mentioned above, one can see that a necessary and
sufficient condition for $\vartheta_p>0$ is 
$\mu_{11}\Hat\nu_1 +\mu_{12} \Hat\nu_2 >0$.
If $\Hat\nu_1<0$ and $\Hat\nu_2>0$, by rewriting the condition as
$\frac{\mu_{11}\Hat\nu_1}{\mu_{12}} +  \Hat\nu_2 >0$,
we see that,
the first term represents the service capacity required for class $1$ at pool $2$ to
be reallocated,  and thus, the sum being positive means that there is an allowance at
pool $2$ for class $1$ to be served. 
Similarly for the case when $\Hat\nu_1>0$ and $\Hat\nu_2<0$.
\end{example}

\begin{example}[The `M' Network]
We obtain the SWSS parameter
\begin{equation*}
\vartheta_p \,=\, \frac{\mu_{22}}{\mu_{22}p_1+ \mu_{12}p_2} 
\Bigl(\mu_{11}\theta_1 + \mu_{12}\theta_2
+ \frac{\mu_{12}\mu_{23}}{\mu_{22}} \theta_3 -\Hat\lambda_1
- \frac{\mu_{12}}{\mu_{22}}\Hat\lambda_2 \Bigr)\,,
\end{equation*}
with 
\begin{equation*}
\kappa_{11}^* \,=\, \theta_1\,, \quad \kappa_{23}^*
\,=\, \theta_3\,, \quad \kappa_{12}^* 
\,=\, \theta_2 - \kappa_{22}^*\,, \quad \kappa_{22}^*
\,=\, \frac{\Hat\lambda_2 + \vartheta_p\, p_2 - \mu_{23} \theta_3}{\mu_{22}}\,.
\end{equation*}
In this case, the $B_1$ matrix is given by $B_1 = \diag (\mu_{12},\mu_{22})$,
and the vector $h$ is given by 
\begin{equation*}
   h \,=\, 
    \begin{pmatrix}
    \Hat{\lambda}_1 - \mu_{11}\xi^*_{11}\Hat{\nu}_1 - \Hat{\mu}_{11} z^*_{11}
    - \mu_{12}\xi^*_{12}\Hat{\nu}_2 - \Hat{\mu}_{12} z^*_{12}\\[3pt]
    \Hat{\lambda}_2 - \mu_{22}\xi^*_{22}\Hat{\nu}_2 - \Hat{\mu}_{22} z^*_{22}
    -\mu_{23}\xi^*_{23}\Hat\nu_3 - \Hat\mu_{23}z^*_{23}
    \end{pmatrix}\,,
\end{equation*}
where $\xi^*_{11} = 1$, $\xi^*_{12} + \xi^*_{22} = 1$, $\xi^*_{23} = 1$.
It is  clear that
$\vartheta_p = -\frac{\langle e,B_1^{-1} h \rangle}{\langle e,B_1^{-1}p\rangle}$.
In the special case, a necessary and sufficient condition for $\vartheta_p>0$ is
\begin{equation*}
\mu_{11}\Hat\nu_1 + \mu_{12}\Hat\nu_2 + \frac{\mu_{12}\mu_{23}}{\mu_{22}}\Hat\nu_3 \,>\,0
\iff \frac{\mu_{11}}{\mu_{12}}\Hat\nu_1 + \Hat\nu_2
+ \frac{\mu_{23}}{\mu_{22}}\Hat\nu_3 \,>\,0\,.
\end{equation*}
This condition also has a very intuitive interpretation.
For instance, if $\Hat\nu_1<0\,, \Hat\nu_3<0$ and $\Hat\nu_2>0$,  
then the first and third terms represent the service capacity required 
for class 1 and class 3 at server pool 2, respectively,
 and the sum being positive means that the safety staffing at pool 2 is 
 sufficient to serve these additional service requirements. 
\end{example}

\begin{example}[The `W' Network]
The SWSS parameter is given by
\begin{equation*}
\vartheta_p \,=\, \frac{1}{\frac{\mu_{21}}{\mu_{11}}p_1 + p_2
+ \frac{\mu_{22}}{\mu_{32}}p_3} \Bigl(\mu_{21}\theta_1
+ \mu_{22}\theta_2 - \frac{\mu_{21}}{\mu_{11}}\Hat\lambda_1 - \Hat\lambda_2
- \frac{\mu_{22}}{\mu_{32}}\Hat\lambda_3\Bigr)\,,
\end{equation*}
with 
\begin{equation*}
\kappa_{11}^* \,=\, \frac{\Hat\lambda_1 + \vartheta_p\, p_1}{\mu_{11}}\,,
\quad \kappa_{21}^*
\,=\, \theta_1 - \kappa_{11}^*\,, \quad \kappa_{32}^*
\,=\, \frac{\Hat\lambda_3 + \vartheta_p\, p_3}{\mu_{32}}\,\quad
\kappa_{22}^* \,=\, \theta_2 - \kappa_{32}^*.
\end{equation*}
In this case, the $B_1$ matrix and $h$ vector are given by
\begin{equation*}
 B_1 \,=\, \begin{pmatrix} \mu_{11} & 0 &0 \\[3pt]
\mu_{22}-\mu_{21}  &  \mu_{22} & 0 \\[3pt] 0 & 0 & \mu_{32} \end{pmatrix}
\qquad\text{and}\qquad 
h \,=\ \begin{pmatrix} 
\Hat\lambda_1 -\mu_{11}\xi^*_{11}\Hat\nu_1 -\Hat\mu_{11}z^*_{11}\\[3pt]
\Hat\lambda_2 - \mu_{21}\xi^*_{21}\Hat\nu_1-\mu_{22}\xi^*_{22}\Hat\nu_2
- \Hat\mu_{21}z^*_{21} -\Hat\mu_{22}z^*_{22}\\[3pt]
\Hat\lambda_3 - \mu_{32}\xi^*_{32}\Hat\nu_3 -\Hat\mu_{32}z^*_{32}
\end{pmatrix}\,,
\end{equation*}
where $\xi^*_{11} + \xi^*_{21} = 1$ and $\xi^*_{22} + \xi^*_{32} = 1$.
A simple calculation confirms that
$\vartheta_p = -\frac{\langle e,B_1^{-1} h \rangle}{\langle e,B_1^{-1}p\rangle}$.
In the special case, a necessary and sufficient condition for $\vartheta_p>0$ is 
$\mu_{21}\Hat\nu_1 + \mu_{22}\Hat\nu_2 >0$.
\end{example}

\section{Transience}\label{S5}
In this part, we show that both the diffusion limit and diffusion-scaled state process
of the $n^{\mathrm th}$ system
are transient
when $\vartheta_p<0$ and $\vartheta^n_p <0$, respectively.
In addition, we show that they cannot be positive recurrent when
$\vartheta_p = 0$ and $\vartheta^n_p = 0$, respectively.
We start with the following important lemma.

\begin{lemma}\label{L2}
The drift in \cref{E-drift2} satisfies
$\inf_{u^s\in\varDelta_s}\,
\bigl(1+\bigl\langle e ,B_1^{-1} B_2 u^s\bigr\rangle\bigr) > 0$,
where $\varDelta_J$ is as in \cref{E-Act}.
\end{lemma}
\proof{Proof.}
This proof is motivated by that of \cref{T3}. Let $\tilde{\vartheta}_p \in \RR_+$,  $\tilde{\kappa}=[\tilde{\kappa}_{ij}]\in\RR_+^{\cG}$ be the unique solution of the following optimization problem
\begin{equation}\label{E-LP''}
\begin{aligned}
\text{maximize}&\quad \tilde{\vartheta}_p \\
\text{subject to}&\quad
0\,\le\,
\sum_{j \in \cJ(i)} \mu_{ij}\tilde\kappa_{ij}-\tilde\vartheta_p\,p_i 
\qquad\forall\,i\in\cI\,,\\
& \sum_{i\in\cI(j)} \tilde\kappa_{ij} \,=\, u^s_j\,,
 \qquad\forall j\in\cJ, u^s \in \varDelta_s\,.
 \end{aligned}\tag{\sf{LP$^{''}$}}
\end{equation}
Note that if we set $\hat{\lambda}_i = 0$, $\theta_j = u^s_j$ in the proof of \cref{T3}, the same conclusion holds as in \cref{EPT3A}, i.e., 
\begin{equation*}
\tilde\vartheta_p e^T B_1^{-1} p  = e^T u^s + e^T B_1^{-1} B_2 u^s = 1 + e^T B_1^{-1} B_2 u^s,
\end{equation*}
where the last equality holds because $u^s \in \varDelta_s$. The proof is concluded by noting that $\tilde\vartheta_p \in \RR_+$ and that we have shown in \cref{T3} that $\bigl(e\transp B_1^{-1}\bigr)_i>0$ for all $i\in\cI$.
\Halmos\endproof

We first show that $\vartheta_p<0$ implies transience for the diffusion limit. 

\begin{theorem}\label{T5}
Suppose that $\vartheta_p<0$.
Then the limiting diffusion $\process{X}$ in \cref{E-diff} is transient under any
stationary Markov control.
In addition, if $\vartheta_p= 0$, then $\{X(t)\}_{t\ge0}$
 cannot be positive recurrent.
\end{theorem}

\proof{Proof.}
In the following, note that 
the function $H(x)$ is a test function and it is chosen such that $\Lg_u H(x) > 0$.
Let	$H(x) \df \tanh\bigl(\beta\langle e, B_1^{-1} x\rangle\bigr)$, with $\beta>0$.
Then
\begin{equation*}
\trace\bigl(\Sigma\Sigma\transp\nabla^2 H(x))\bigr) \,=\,
\beta^2\tanh''\bigl(\beta\langle e, B_1^{-1} x\rangle\bigr)
\babs{\Sigma\transp B_1^{-1}e}^2\,,
\end{equation*}
where we recall $\Sigma \df \diag\bigl(\sqrt{2 \lambda_1}, \dotsc, \sqrt{2 \lambda_I}\bigr)$.
We have
\begin{equation}\label{PT5A}
\begin{aligned}
\Lg_u H(x) &\,=\, \frac{1}{2}\trace\bigl(\Sigma\Sigma\transp\nabla^2 H(x)\bigr)
+ \bigl\langle b(x,u),\nabla H(x)\bigr\rangle \\
& \,=\, - \beta^2\frac{\tanh\bigl(\beta\langle e, B_1^{-1} x\rangle\bigr)}
{\cosh^2\bigl(\beta\langle e, B_1^{-1} x\rangle\bigr)}
\abs{\Sigma\transp B_1^{-1} e}^2\\
&\mspace{80mu}+ \frac{\beta}{\cosh^2\bigl(\beta\langle e, B_1^{-1} x
\rangle\bigr)}\biggl(\bigl\langle e, B_1^{-1} h\bigr\rangle 
+ \langle e,x\rangle^-
\Bigl(1 + \bigl\langle e , B_1^{-1}B_2 u^s\bigr\rangle\Bigr)\biggr)\,.
\end{aligned}
\end{equation}
Thus, for
$0<\beta < \langle e, B_1^{-1} h\rangle\, \abs{\Sigma\transp B_1^{-1} e}^{-2}$,  
we obtain $\Lg_u H(x) >0$ by \cref{L2} and using \cref{T3,T4} to conclude that
$\langle e,B_1^{-1} h \rangle >0$. Therefore, $\{H\bigl(X_t\bigr)\}_{t\ge0}$
is a bounded submartingale,
so it converges almost surely.
Since $X$ is irreducible, it can be either recurrent or transient. 
If it is recurrent, then $H$
should be constant a.e.\ in $\RR^I$, which is not the case.
Thus $X$ is transient.

We now turn to the case where $\vartheta_p = 0$.
Suppose that the process $\{X(t)\}_{t\ge0}$ (under some stationary
Markov control) has an invariant probability measure $\uppi(\D{x})$.
It is well known that $\uppi$ must have a positive density.
Let $g_1(x)$ and $g_2(x)$ denote respectively the first and the second terms on
the right hand side of \cref{PT5A}.
Applying It\^{o}'s formula to \cref{PT5A}, we obtain
\begin{equation}\label{PT5B}
\Exp^\uppi\bigl[H(X_{t\wedge \tau_r})\bigr] - H(x) \,=\,
\sum_{i=1,2}\Exp^\uppi\Biggl[\int_0^{t\wedge \tau_r} g_i(X_s)\D s\Biggr]\,,
\end{equation}
where $\tau_r$ denotes the first exit time from the ball $B_r$ of radius $r$
centered at $0$.
Note that $g_1(x)$ is bounded and $g_2(x)$ is non-negative.
Thus using dominated and monotone convergence, we can take limits in \cref{PT5B}
as $r\to \infty$ for the terms on the right side to obtain
\begin{equation*}
\int_{\Rm} H(x)\uppi(\D{x}) - H(x) \,=\,
t \sum_{i=1,2} \int_{\Rm} g_i(x)\uppi(\D{x}), \qquad t\ge0\,.
\end{equation*}
Since $H(x)$ is bounded, we can divide both sides by $t$ and $\beta$ and take the
limit as $t\to\infty$ to get
\begin{equation}\label{PT5C}
\int_{\Rm} \beta^{-1}g_1(x) \uppi(\D{x})
+ \int_{\Rm} \beta^{-1}g_2(x) \uppi(\D{x}) \,=\, 0\,.
\end{equation}
Since $\beta^{-1} g_1(x)$ tends to 0 uniformly in $x$ as $\beta\searrow 0$,
the first term on the left hand side of \cref{PT5C} vanishes as
$\beta\searrow0$.
However, since $\beta^{-1}g_2(x)$ is bounded away from $0$ on the open set
$\{x\in\Rm \colon \langle e,x \rangle ^- >1\}$,
this contradicts the fact that $\uppi(\D{x})$ has full support.
\Halmos\endproof

The proof of the following corollary is analogous to that of \cref{L2}.

\begin{corollary}\label{C2}
The drift in \cref{E-bn} satisfies
$\inf_{u^s\in\varDelta_s}\,
\bigl(1+\bigl\langle e ,(B_1^n)^{-1} B_2^n u^s\bigr\rangle\bigr) > 0$.
\end{corollary}

\begin{theorem}\label{T6}
Suppose that $\vartheta_p^n <0$.
Then the diffusion-scaled state process $\{\Breve X^n(t)\}_{t\ge0}$ of the $n^{\mathrm{th}}$
system is transient under any stationary Markov scheduling policy.
In addition, if $\vartheta_p^n = 0,$ the process $\{\Breve X^n(t)\}_{t\ge0}$  cannot be positive recurrent.
\end{theorem}

\proof{Proof.}
The proof mimics that of \cref{T5}. 
We apply the function $H(x) = \tanh\bigl(\beta\langle e, (B_1^n)^{-1} x\rangle\bigr)$
to the operator
$\Breve\Lg^n_{z}$ in \cref{E-brLgn}, and use the identity
\begin{equation*}
H\Bigl(x\pm \tfrac{1}{\sqrt n} e_i\Bigr) - H(x) \mp
\frac{1}{\sqrt n}\partial_{x_i} H(x)
\,=\, \frac{1}{n}\int_0^1 (1-t)\,
\partial_{x_ix_i} H\Bigl(x\pm \tfrac{t}{\sqrt n} e_i\Bigr) \,\D{t}
\end{equation*}
to express the first and second order incremental quotients,
together with \cref{E-bn} which implies that
\begin{multline*}
\bigl\langle b^n(\Breve{x},\Breve{z}),\nabla H(\Breve{x})\bigr\rangle
\,=\,
\frac{\beta}{\cosh^2\bigl(\beta\langle e, (B_1^n)^{-1} \Breve{x}
\rangle\bigr)}\biggl(\bigl\langle e, (B_1^n)^{-1}h^n\bigr\rangle \\
+ \bigl(\Breve{\zeta}^{n}(\Breve{x},\Breve{z})+\langle e,\Breve x\rangle^-\bigr)
\Bigl(1 + \bigl\langle e , (B_1^n)^{-1}B_2^n u^s\bigr\rangle\Bigr)\biggr)\,.
\end{multline*}
The rest follows exactly as in the proof of \cref{T5} using \cref{C2}.
\Halmos\endproof

\section{Stability}\label{S6}
We start with the following important lemma which is essential in
proving the stabilization results.

\begin{lemma}\label{L5}
Suppose that the solution $\vartheta_p$ of \cref{E-LP'}
is positive.
Then, there exist a collection $\{\widetilde{N}_{ij}^n\in\NN,(i,j)\in\cE, n\in\NN\}$,
and a positive constant $C_0$ satisfying 
\begin{align}
 \lambda_i^n &\,= \, \sum_{j \in \cJ(i)} \mu^n_{ij}  \widetilde{N}_{ij}^n 
- \vartheta_p\,  p_i \sqrt{n} + \sorder\bigl(\sqrt n\bigr)\,.
\qquad\forall\,i\in\cI\,,\label{EL5A} \\[5pt]
\babs{N_{ij}^n-\widetilde{N}_{ij}^n} &\,\le\,
C_0\sqrt{n}\qquad\forall (i,j)\in\cE\,,\label{EL5B}
\end{align}
with $\widetilde{N}_{ij}^n$ as in \cref{E-Nij}, and
\begin{equation}\label{EL5C}
\sum_{i\in \cI(j)}\widetilde{N}_{ij}^n \,=\, N_j^n\qquad\forall\,j\in\cJ\,,
\end{equation}
for all sufficiently large $n\in\NN$. 
\end{lemma}

\proof{Proof.}
Let $\{\kappa_{ij}\}_{(i,j)\in\cE}$ be a solution of the optimization problem
in \cref{E-LP'}.
There is flexibility in selecting such
a set $\{\widetilde{N}_{ij}^n\}_{(i,j)\in\cE}$.
For example, first select some arbitrary
element $\Hat\imath_j\in\cI(j)$ each $j\in\cJ$.
It is clear that we can select a set of numbers
$\widetilde{N}_{ij}^n$, $(i,j)\in\cE$, satisfying \cref{EL5C},
which also satisfy
\begin{equation} \label{PL5A}
\biggl \lfloor n z^*_{ij} + \sqrt{n} \biggl(\kappa_{ij}
- \frac{\Hat\mu_{ij}}{\mu_{ij}} z^*_{ij} \biggr) \biggr\rfloor
\,\le\, \widetilde N_{ij}^n
\,\le\, \biggl\lceil n z^*_{ij}+ \sqrt{n} \biggl(\kappa_{ij}
- \frac{\Hat\mu_{ij}}{\mu_{ij}} z^*_{ij} \biggr)\biggr\rceil
\qquad\forall\,i\in\cI\setminus\{\Hat\imath_j\}\,,\ \forall\,j\in\cJ\,.
\end{equation}
Then \cref{EL5B} holds by construction.
Using \cref{E-LP',PL5A}
in combination with $\sum_{j \in \cJ} \mu_{ij} z^*_{ij} = \lambda_i$
and the convergence of parameters in \cref{E-Par01,E-Par02}, it is easy to see
that \cref{EL5A} holds.
\Halmos\endproof

Let $\widetilde{N}^n_i \df \sum_{j\in \cJ(i)}\widetilde{N}_{ij}^n$
for $i\in\cI$, and define
\begin{equation}\label{E-tVar}
\Tilde{x}^{n}_i(x) \,\df\, \frac{1}{\sqrt{n}} \bigl(x -\widetilde{N}^n_i\bigr)\,.
\end{equation}
Recall that the matrices $B_1$ and $B_2$ in \cref{E-drift2}
are independent of the choice of centering.
Thus, employing the same approach
as in \cref{S4}, one can easily show that
the process $\widetilde{X}^n = \Tilde{x}^n(X^n)$ converges to the limit $X$ described
in \cref{E-diff} with $h_i = - \vartheta_p \,p_i$ for all $i\in\cI$.
It also follows from \cref{L5} that the expression in \cref{E-hn} gets replaced by
\begin{equation*}
h_i^n \,=\, \frac{1}{\sqrt{n}}
\biggl(\lambda_i^n - \sum_{j \in \cJ(i)} \mu^n_{ij}  \widetilde{N}_{ij}^n\biggr)
\,=\, -\vartheta_p \,p_i
+ \frac{\sorder\bigl(\sqrt n\bigr)}{\sqrt{n}}\,.
\end{equation*}
We emphasize here that
if $\vartheta_p>0$, then
under this rebalancing mechanism, 
 the aggregate steady-state capacity
$\sum_{j \in \cJ(i)} \mu^n_{ij}  \widetilde{N}_{ij}^n $ 
provides sufficient safety staffing for each class $i$ by \cref{EL5A}.

\subsection{Stabilizing the Limiting Diffusion}\label{S6.1}

We define the following class of Markov controls
$\Bar{v}=(\Bar{v}^c,\Bar{v}^s)$ for the diffusion.
Let $(\Hat\imath,\Hat\jmath)\in\cE$ be given. Define
\begin{equation}\label{E-barv}
\Bar{v}^c(x) \,=\, e_{\Hat\imath}\quad\text{and\ } \Bar{v}^s(x) \,=\, e_{\Hat\jmath}\,,
\quad x\in\Rm\,.
\end{equation}
The control $\Bar{v}$ can be interpreted as follows:
All the jobs in queue, if any, are in class $\Hat\imath$, and any idle servers
are in pool $\Hat\jmath$.
This interpretation of course applies to the limiting diffusion,
but note that such a control is admissible for the $n^{\mathrm th}$ system
in the set \cref{E-jwc}.
Thus, if this control is applied to the $n^{\mathrm th}$ system,
then it is clear that the drift in \cref{E-bn} converges, as $n\to\infty$, to
the drift
the limiting diffusion controlled under $\Bar{v}$.

Consider the diffusion limit in \cref{E-drift2}. Recall that $h_i = -\vartheta_p p_i$ under the new centering used in \cref{E-tVar}.
Choosing an ordered
$\mathfrak{D} =\bigl(\alpha, (\beta)_{-\Hat\jmath}\bigr)$,
with $\alpha_{\Hat\imath}$ being the last element of $\alpha$,
and applying \cref{P2},
we see that the drift of the diffusion takes the form
\begin{equation*}
\Bar{b}(x)\,\df\, b\bigl(x,\Bar{v}(x)\bigr) \,=\,
\begin{cases}
-\vartheta_p\,p
- B_1 (\Id -  e_I e\transp)x\,, &\text{if\ } \langle e, x\rangle\ge0\,,\\[3pt]
-\vartheta_p\,p - B_1 x\,, &\text{if\ } \langle e, x\rangle < 0\,.
\end{cases}
\end{equation*}
Note that the term $B_2\Bar{v}^s(x) \langle e,x \rangle^-$ does not appear
in the representation of $\Bar{b}(x)$ above when $ \langle e, x\rangle < 0$.
This is because the last column of $B_2$ is identically zero as noted in \cref{P2}.

Note that $B_1$ and $B_1(\Id -  e_I  e\transp)$ 
are both lower diagonal matrices, where $B_1(\Id -  e_I  e\transp)$
has all positive diagonal elements except for the $I^{\rm th}$ one which equals zero.
Therefore, $B_1^2(\Id -  e_I  e\transp)$ has no real negative eigenvalues
and a simple zero eigenvalue.
Thus, by \citet[Proposition 3]{DG13}, there exist a positive definite matrix
$S\in\RR^{I\times I}$, and a constant $\kappa_\circ>0$, such that
\begin{equation*} 
S B_1 + B_1\transp S \,>\, 2\kappa_\circ\Id\,, \quad
\text{and}\quad  \Phi\,\df\, S B_1(\Id -  e_I e\transp)
+ (\Id- e  e_I\transp) B_1\transp S \,\ge\, 0\,.
\end{equation*}
Define
$\eta\df \vartheta_p \langle p, S  e_I\rangle$.
Note that, since $S$ is a positive
definite matrix, there exists a positive vector $p\in\varDelta_I$ such that $\eta>0$.

Let
$\norm{x}_S \df \langle x,Sx\rangle^{\nicefrac{1}{2}}$, and define 
\begin{equation*}
\Lyap_{\epsilon,S}(x) \,\df\,
\exp\Bigl( \epsilon \norm{x}_S^2\,
\bigl(1+\norm{x}_S^2\bigr)^{-\nicefrac{1}{2}}\Bigr)\,,\qquad x\in\RR^I\,.
\end{equation*}

Recall \cref{E-Lg}. We have the following theorem.

\begin{theorem}\label{T7}
Fix $p\in\varDelta_I$ such that $\eta>0$ and assume that $\vartheta_p >0$. Let $\Bar{v}$ be as in
\cref{E-barv}. Then, there exist
a positive definite matrix $S\in\RR^{I\times I}$
and positive constants $\epsilon$ and $\kappa_i$, $i=0,1$, such that
\begin{equation}\label{ET7A}
\Lg_{\Bar{v}}\,\Lyap_{\epsilon,S}(x) \,\le\,
\kappa_0 - \kappa_1\,\Lyap_{\epsilon, S}(x)\qquad\forall\,x\in\RR^{I}\,.
\end{equation}
The process $\process{X}$ is exponentially ergodic and admits a unique
invariant probability measure $\uppi_{\Bar{v}}$ under $\Bar{v}$ satisfying
\begin{equation}\label{ET7B}
 \bnorm{P^{\Bar{v}}_t (x, \cdot)
- \uppi_{\Bar{v}}(\cdot)}_{\Lyap_{\epsilon,S}} \,\le\,
C_\gamma \Lyap_{\epsilon,S} (x) e^{-\gamma t} \,,
\quad x \in \RR^I\,, \quad \forall t\,\ge\,0\,,
\end{equation}
where $P^{\Bar{v}}_t (x, dy)$ denotes the transition probability of $\process{X}$
under $\Bar{v}$.
\end{theorem}
Before proving the theorem, we would like to recall \cref{R1}. Note that introducing the vector $p$ and fixing it throughout the proof of \cref{T7} is critical. This is mainly the reason behind using $p$ when defining the SWSS $\vartheta_p$ instead of just using the numerator of the expression in \cref{ET2A}.
\proof{Proof.}
%
Define $\delta\df\frac{1}{4}\kappa_\circ\,\abs{SB_1  e_I}^{-1}$ and recall that
$\eta = \vartheta_p \langle p, S  e_I\rangle$.
Let $\varphi(x) \df \langle x,Sx\rangle$, and
$\cK_{\delta} \df \bigl\{x\in\RR^I\,\colon \langle e,x\rangle
> \delta \abs{x}\bigr\}$.
For $x\in\cK_\delta^c$,  we obtain
\begin{align*}
\bigl\langle \Bar{b}(x), \nabla \varphi(x)\bigr\rangle 
&\,=\, - 2 \vartheta_p \langle p, Sx\rangle - \langle x,(SB_1+B_1\transp S)x\rangle
+2\langle x, SB_1  e_I\rangle \langle e, x\rangle^+\\
&\,\le\, 2 \vartheta_p\abs{S p}\abs{x} -2\kappa_\circ \abs{x}^2
+2 \delta \abs{SB_1  e_I}\abs{x}^2\,.
\end{align*}
Thus, by the definition of $\delta$,
and with $\Bar{\kappa}\df 2 \kappa_\circ^{-1}\bigl(\vartheta_p\abs{S p}\bigr)^2$,
we have
\begin{equation*}
\bigl\langle \Bar{b}(x), \nabla \varphi(x)\bigr\rangle
\,\le\, \Bar{\kappa} - \kappa_\circ \,\abs{x}^2\,,\qquad \forall\,x\in\cK_\delta^c\,.
\end{equation*}

Next, suppose that $x\in\cK_\delta$.
We have 
\begin{equation}\label{PT7A}
\bigl\langle \Bar{b}(x), \nabla \varphi(x)\bigr\rangle 
\,=\, -2 \vartheta_p \langle p, Sx\rangle - \bigl\langle x, \Phi x\bigr\rangle \,. 
\end{equation}
Decompose $x = x_{(-I)} + x_I  e_I$ into the orthogonal components
$x_{(-I)}$ and $x_I  e_I$. 
Then
\begin{equation}\label{PT7B}
\begin{aligned}
\vartheta_p  \langle p, Sx\rangle &\,=\,  \vartheta_p 
\langle p, S x_{(-I)}\rangle +      \eta x_I \,, \\
&\,=\,   \vartheta_p   \langle p, S x_{(-I)}\rangle
+ \eta \Bigl( \langle e,x\rangle - \langle e,  x_{(-I)} \rangle \Bigr)  \, \\
&\,\ge\,   \vartheta_p   \langle p, S x_{(-I)}\rangle
- \eta \langle e,  x_{(-I)} \rangle  + \eta \delta \abs{x} \,,
\end{aligned}
\end{equation}
and
\begin{equation}\label{PT7C}
\begin{aligned}
\bigl\langle x, \Phi x\bigr\rangle  &\,=\, \bigl\langle x_{(-I)} + x_I  e_I,
\Phi \bigl( x_{(-I)} + x_I  e_I\bigr)  \bigr\rangle  \\
&\,=\,  x\transp_{(-I)} \Phi x_{(-I)}  + 2  x_I  e_I\transp
\Phi  x_{(-I)} + x_I^2  e_I\transp \Phi  e_I  \\
&\,=\,  x\transp_{(-I)} \Phi x_{(-I)}\,,
\end{aligned}
\end{equation}
where the last equality uses the fact that $ e_I\transp   \Phi  e_I  = 0$
which implies that $e_I\transp \Phi = 0$ since $\Phi$ is a positive
semi-definite matrix,
and which in its turn implies that $e_I\transp S B_1(\Id -  e_I e\transp) = 0$.
This implies that
\begin{equation*}
e_I\transp S B_1 = \bigl(e_I\transp S B_1 e_I\bigr)\,e\transp \iff
e_I\transp S \,=\, \bigl(e_I\transp S B_1 e_I\bigr)\, e\transp B_1^{-1}\,.
\end{equation*}
Thus, arguing as in the derivation of (5.18)--(5.19) in \citet{DG13},
we conclude that
\begin{equation}\label{PT7D}
\begin{aligned}
x_{(-I)}\transp \Phi x_{(-I)} &\,=\, x_{(-I)}\transp
\bigl(\Id - e e_I\transp\bigr)B_1\transp
\Bigl(S B_1^{-1} + (B_1^{-1})\transp S\Bigr) 
B_1 \bigl(\Id - e_I e\transp\bigr) x_{(-I)} \\
&\,\ge\, c \abs{x_{(-I)}}^2\,
\end{aligned}
\end{equation}
for some positive constant $c$,
where in the last inequality we used the fact that the zero eigenvalue
of $\Id - e_I e\transp$
is simple, and the corresponding eigenvector is $e_I$.
Combining \cref{PT7A,PT7B,PT7C,PT7D}, we obtain 
\begin{equation*}
\bigl\langle \Bar{b}(x), \nabla \varphi(x)\bigr\rangle
\,\le\, \Hat\kappa_0 - \delta\eta\abs{x}
\qquad \forall\,x\,\in\cK_\delta\,,
\end{equation*}
for some constant $\Hat\kappa_0>0$.

Next, if we let 
$\phi_S(x) \df \frac{2+\langle x,Sx\rangle}{(1+\langle x,Sx \rangle)^{\nicefrac{3}{2}}}$,
then a straightforward calculation shows that
\begin{equation*}
\nabla \Lyap_{\epsilon,S}(x) \,=\,
\frac{1}{2} \epsilon\,\Lyap_{\epsilon,S}(x)\,\phi_S(x)\,\nabla \varphi(x)\,,
\end{equation*}
and
\begin{align*}
\nabla^2 \Lyap_{\epsilon,S}(x) &\,=\,
\epsilon^2\,\Lyap_{\epsilon,S}(x)\,\phi^2_S(x)\, Sx\transp xS
+ \epsilon \Lyap_{\epsilon,S}(x) \biggl[\phi_S(x) S
+ \frac{Sx\transp xS}{(1+x\transp Sx)^\frac{5}{2}} \bigl(-4 -\norm{x}_S^2\bigr) \biggr]\\
&\,\le\, \epsilon^2\,\Lyap_{\epsilon,S}(x)\,\phi^2_S(x)\, Sx\transp xS
+ \epsilon \Lyap_{\epsilon,S}(x)\,\phi_S(x) S\,.
\end{align*}
Therefore, if we choose $\epsilon>0$ small enough,
then for some positive constants $\kappa_0$ and $\kappa_1$
we obtain
\begin{equation*}
\begin{aligned}
\Lg_{\Bar{v}}\,\Lyap_{\epsilon,S}(x)&\,=\,
\frac{1}{2}\trace \bigl(\Sigma\Sigma\transp \nabla^2 \Lyap_{\epsilon,S}(x) \bigr)
+ \bigl \langle \Bar{b}(x), \nabla \Lyap_{\epsilon,S}(x)\bigr \rangle
\,\le\, \kappa_0 - \kappa_1 \Lyap_{\epsilon,S}(x)\,,
\end{aligned}
\end{equation*}
which establishes \cref{ET7A}.
It is well known that this drift inequality implies \cref{ET7B}
(see \citet[Theorems~4.3 and 6.1]{MT-III-93}, or \citet[Theorem~5.2]{DMT-95}).
This completes the proof.
\Halmos\endproof

{\color{dred}
\begin{remark}\label{R7}
Under any control $v\in\Usm$ which renders the diffusion limit
positive recurrent with invariant probability measure $\uppi_v$
we have the following identity:
\begin{equation*}
\int_{\Rm} \Bigl(1 + \bigl\langle e,B_1^{-1}B_2v^s(x)\bigr\rangle\Bigr)
\langle e,x \rangle ^-\, \uppi_v(\D{x})
\,=\, \langle e B^{-1}_1, p \bigr\rangle\,\vartheta_p\qquad\forall\, p\in\varDelta_I\,.
\end{equation*}
This extends \citet[Theorem 3.1]{HAP19} to arbitrary tree topologies.
In particular, for the control $\Bar{v}$ in \cref{E-barv}
we obtain 
\begin{equation*}
\int_{\Rm} \langle e,x \rangle ^-\, \uppi_{\Bar{v}}(\D{x})
\,=\, \langle e B^{-1}_1, p \bigr\rangle\,\vartheta_p\,.
\end{equation*} 
This can be interpreted as follows.
The average number of idle servers under the control $\Bar{v}$
equals $\langle e B^{-1}_1, p \bigr\rangle\,\vartheta_p$.

We also mention, parenthetically, that
$\int_\Rm \Lyap_{\epsilon,S}(x)\,\mu(\D{x}) \le \frac{\kappa_0}{\kappa_1}$
by \cref{ET7A}.
This means that the invariant probability measure
$\uppi_{\Bar{v}}$ has exponential tail.
\end{remark}
}

\begin{remark}\label{R8}
The use of $X$ to denote the diffusion limit of $\Breve{X}^n$ and $\tilde{X}^n$ is just an abuse of notation. As mentioned in \cref{R5}, the limiting processes using different centering terms
only differ in the drift by a constant,
and therefore they are equivalent as far as their ergodic properties
are concerned. Here is a more detailed explanation: 

Let $\Breve{X}$ and $\tilde{X}$ denote the limiting processes of $\Breve{X}^n$ and $\tilde{X}^n$ respectively. Hence we have the following SDEs
\begin{equation*}
\begin{aligned}
d \breve{X}_t \,&=\, \breve b( \breve{X}_t, U_t)\,\D{t} + \Sigma \, \D W_t\,,\\
d \tilde{X}_t \,&=\, \tilde b( \tilde{X}_t, U_t)\,\D{t} + \Sigma \, \D W_t\,,
\end{aligned}
\end{equation*}
where
\begin{equation*}
\begin{aligned}
\breve b(\breve{x}, u) \,&=\, h - B_1 (\breve{x} -  \langle e, \breve{x} \rangle^{+} u^c)
+ \langle e, \breve{x} \rangle^{-} B_2 u^s\,,\\
\tilde b(\tilde{x}, u) \,&=\, -\vartheta_p p - B_1 (\tilde{x} -  \langle e, \tilde{x} \rangle^{+} u^c)
+ \langle e, \tilde{x} \rangle^{-} B_2 u^s\,.
\end{aligned}
\end{equation*}
Define $\zeta = B_1^{-1} h + \vartheta_p B_1^{-1}p$. Using \cref{T4}, we have $\langle e, \zeta \rangle = \langle e, B_1^{-1} h \rangle + \vartheta_p \langle e, B_1^{-1}p \rangle = 0$. One can then check that $\breve{X} = \tilde{X} + \zeta$ and hence the ergodic properties of the limiting processes are the same.

\end{remark}

\begin{remark}
We note that when $\vartheta_p>0$, the class
of stabilizing controls might be much richer.
Indeed, it has been shown in \citet{HAP19} that if
$\varrho = -\langle e,B_1^{-1} h \rangle>0$, where $h$ is given by \cref{E-h},
the diffusion limit of networks with a single non-leaf server pool and those whose
service rates are dictated by the class type are uniformly exponentially ergodic
under any stationary Markov control.
In addition, the prelimit diffusion-scaled processes are uniformly exponentially
ergodic over a class of policies which is referred to as system-wide work-conserving
in \citet{HAP19}.
Using the equivalence relation between $\varrho$ and $\vartheta_p$ in \cref{S4},
these conclusions hold for these networks when $\vartheta_p>0$.
\end{remark}

\subsection{Stabilizing the Diffusion-Scaled Processes} \label{S6.2}

Recall \cref{E-tVar,E-Lgn}. 
The generator $\widetilde\Lg^n_{z}$ of the diffusion-scaled
state process $\widetilde{X}^n$  satisfies 
\begin{equation} \label{E-tLg} 
\widetilde\Lg^n_{z} f(y)\bigr|_{y=\Tilde{x}^n(x)}
\,=\, \Lg^n_{z} \bigl(f\circ \Tilde{x}^n\bigr)(x)\,.
\end{equation}

A family of scheduling policies, referred to as
\emph{balanced saturation policies} (BSPs), is introduced in \citet{AP19}.
When there is at least one class with positive abandonment rate, exponential ergodicity
is shown under the BSPs (see Proposition 5.1 therein). 
The proof of this result relies on the system having a positive
abandonment rate in some class, and cannot be applied directly here.
Provided that $\vartheta_p>0$, we we show in \cref{T8}
that the diffusion-scaled
processes controlled by a BSP are exponentially ergodic for networks without abandonment,
and the corresponding stationary distributions are tight.
 Recall the definition of the BSPs. 

\begin{definition}
Let $\{\widetilde{N}_{ij}^n\}$ be as in \cref{L5},
and recall that $\widetilde{N}_i^n =  \sum_{j\in \cJ(i)} \widetilde{N}_{ij}^n$
for $i\in\cI$.
Let $\sZ^n$ denote the class of work-conserving Markov policies
$z$ satisfying
\begin{equation}\label{ED6.1}
\begin{aligned}
z_{ij}(x) &\,\le\, \widetilde{N}_{ij}^n \quad\forall\, i\sim j\,,
\quad\text{and}~ \sum_{j\in\cJ(i)} z_{ij}(x)=x_i\,,
\qquad \text{if~} x_i \le \widetilde{N}_i^n\,,\\
z_{ij}(x) &\,\ge\, \widetilde{N}_{ij}^n \quad \forall\, i\sim j\,,
\qquad \text{if~} x_i > \widetilde{N}_i^n\,.
\end{aligned}
\end{equation} 
\end{definition}

It is rather simple to verify that the class $\sZ^n$ is nonempty.
For example, a policy in $\sZ^n$ can be determined in two steps.
In the first step, if $x_i>\widetilde{N}_i^n$, then we set
$z_{ij}(x)=\widetilde{N}_{ij}^n$ for all $j\in\cJ(i)$; otherwise
we determine $z_{ij}(x)$ in any arbitrary manner that satisfies
\cref{ED6.1}.
In the second step, we fill in the pools in any arbitrary manner
that enforces work conservation.
{\color{dred}The following examples of BSPs are adapted from \citet[Definition~3.1]{AP18}
and  \citet[Section~5]{AP19}. 

\begin{example}
We provide explicit definitions of a BSP policy for the `N' and `M' networks.
For the `N' network, this is given by
\begin{align*}
z_{11}(x) & \,=\, x_1\wedge N^n_1 \\
z_{12} (x) &\,=\, 
\begin{cases}
(x_1 - N^n_1)^+ \wedge \widetilde{N}^n_{12} \qquad & \text{if } x_2 \ge \widetilde{N}^n_{22}\\
(x_1 - N^n_1)^+ \wedge (N_2^n - x_2) \qquad &\text{otherwise,}
\end{cases} \\
z_{22}(x) & \,=\, 
\begin{cases}
x_2 \wedge \widetilde{N}^n_{22} \qquad &\text{if } x_1 \ge N_1^n + \widetilde{N}^n_{12}\\
x_2\wedge\bigl(N_2^n - (x_1 - N_1^n)^+\bigr)\qquad & \text{otherwise.}
\end{cases}
\end{align*}
Note that we have used $N^n_1 = \widetilde{N}^n_{11}$.

For the `M' network, a BSP policy is given by:
\begin{align*}
z_{11}(x) & \,=\, x_1\wedge N^n_1 \\
z_{12} (x) &\,=\, 
\begin{cases}
(x_1 - N^n_1)^+ \wedge \widetilde{N}^n_{12} \qquad & \text{if } x_2 \ge \widetilde{N}^n_{2}\\
(x_1 - N^n_1)^+ \wedge (x_2 - N_3^n)^+ \qquad & \text{otherwise,}
\end{cases} \\
z_{22}(x) & \,=\, 
\begin{cases}
(x_2 - N_3^n)^+ \wedge \widetilde{N}^n_{22} \qquad &\text{if } x_1 \ge \widetilde{N}_1^n\\
(x_2 - N_3^n)^+ \wedge(x_1 - N_1^n)^+\qquad & \text{otherwise.}
\end{cases}\\
z_{23}(x) &\,=\, x_2 \wedge N_3^n.
\end{align*}
\end{example}
}

\begin{definition}\label{D6.2} 
For $\epsilon >0$, we define
 \begin{equation*}
\Lyap_{\epsilon}(x) \,\df\,
\exp\Bigl( \epsilon \abs{x}^2
\bigl(1+\abs{x}^2\bigr)^{-\nicefrac{1}{2}}\Bigr)\,,\qquad x\in\RR^I\,,
\end{equation*}
and let $\widetilde\Lyap_{\epsilon}(x) \df\Lyap_{\epsilon}\bigl(\Tilde{x}^n(x)\bigr)$.
\end{definition}

\begin{theorem}\label{T8}
{\color{dred}If $\vartheta_p > 0$, then
there exist $\epsilon>0$, $n_0\in\NN$,
and positive constants $C_0$ and $C_1$ such that
\begin{equation*}
\widetilde\Lg^n_z\, \widetilde\Lyap_{\epsilon}(x) \,\le\,
C_0  - C_1\, \widetilde\Lyap_{\epsilon}(x)\qquad
\text{for all\ } z\in\sZ^n\,,\  x\in\RR^I\,, \text{\ and\ } n\ge n_0\,,
\end{equation*}
with $\widetilde\Lg_n^z$ and $\widetilde\Lyap_{\epsilon}$
as in \cref{E-tLg} and \cref{D6.2}, respectively.}
In particular,
the process $\widetilde{X}^n$ is exponentially ergodic and admits a unique
invariant probability measure $\widetilde\uppi^n$ satisfying
\begin{equation*}
\lim_{t\to\infty} e^{\kappa t} \bnorm{P^n_t (x, \cdot)
- \widetilde\uppi^n(\cdot)}_{\mathrm{TV}} \,=\, 0\,, \quad x \in \RR^I\,,
\end{equation*}
for some $\kappa>0$,
where $P^n_t (x, \cdot)$ denotes the transition probability of $\widetilde{X}^n$. 
\end{theorem}

\proof{Proof.}
Using the identity
\begin{equation*}
f(x\pm e_i) - f(x) \mp \partial_i f(x)
\,=\, \int_0^1 (1-t)\, \partial_{ii} f(x\pm t e_i) \,\D{t}\,,
\end{equation*}
we obtain
\begin{equation}\label{PT8A}
\begin{aligned}
\Babs{\widetilde\Lyap_{\epsilon}(x\pm e_i)
- \widetilde\Lyap_{\epsilon}(x) \mp 
\tfrac{\epsilon}{\sqrt{n}} \Tilde{x}^n_i \phi(\Tilde{x}^n)\, &
\widetilde\Lyap_{\epsilon}(x)} \,\le\,
\tfrac{1}{n}\epsilon^2\,\Tilde\kappa_1\,\widetilde\Lyap_{\epsilon}(x)
\end{aligned}
\end{equation}
for some constant $\Tilde\kappa_1>0$, and all $\epsilon\in(0,1)$, with
$\phi(x) \df \frac{2+\abs{x}^2}{(1+\abs{x}^2)^{\nicefrac{3}{2}}}$.

Fix $n\in\NN$.
Using \cref{PT8A}, we obtain
\begin{align*}
\Lg^n_z\, \widetilde\Lyap_{\epsilon}(x) &\,\le\, \epsilon
\sum_{i\in\cI} \biggl[\lambda_i^n
\Bigl(\tfrac{1}{\sqrt{n}} \Tilde{x}^n_i \phi(\Tilde{x}^n)
+\tfrac{1}{n}\epsilon\,\Tilde\kappa_1\Bigr) + \sum_{j\in \cJ(i)}\mu_{ij}^n z_{ij}\,
\Bigl(-\tfrac{1}{\sqrt{n}} \Tilde{x}^n_i \phi(\Tilde{x}^n)
+\tfrac{1}{n}\epsilon\,\Tilde\kappa_1\Bigr)\biggr]
\widetilde\Lyap_{\epsilon}(x)
\nonumber \\[5pt]
&\,=\, \epsilon\,\widetilde\Lyap_{\epsilon}(x)\,
\sum_{i\in\cI} \Bigl(\tfrac{1}{\sqrt{n}}\,
\phi(\Tilde{x}^n)\,F_{n,i}^{(1)}(x)
+\tfrac{1}{n}\epsilon\,\Tilde\kappa_1\,F_{n,i}^{(2)}(x)\Bigr)\,,
\end{align*}
where
\begin{equation}\label{PT8C}
F_{n,i}^{(1)}(x) \,\df\, 
\Tilde{x}^n_i
\biggl(\lambda_i^n - \sum_{j \in \cJ(i)} \mu_{ij}^n z_{ij}
 \biggr)\,,\quad\text{and}\quad
F_{n,i}^{(2)}(x) \,\df\,  \lambda_i^n
+ \sum_{j \in \cJ(i)} \mu_{ij}^n z_{ij} \,.
\end{equation}
By \cref{EL5A}, there exists  some constant
$\Tilde\kappa_2$ such that for all $n\in\NN$,
\begin{equation}\label{PT8D}
\frac{1}{n} \Bigl(\lambda_i^n+ \sum_{j \in \cJ(i)} \mu_{ij}^n 
\,\widetilde{N}_i^n\Bigr) \,\le\, \Tilde\kappa_2\qquad\forall\,i\in\cI\,. 
\end{equation}
Since  $z_{ij} \le x_i$ for all  $(i,j)\in\cE$, by \cref{PT8C,PT8D}, we obtain
\begin{equation*}
\begin{aligned}
F_{n,i}^{(2)}(x) &\,\le\,  \lambda_i^n + \Bigl( \sum_{j \in \cJ(i)}
\mu_{ij}^n \Bigr) x_i   \\[5pt]
&\,=\,   \lambda_i^n + \Bigl( \sum_{j \in \cJ(i)} \mu_{ij}^n  \Bigr) 
(\widetilde{N}_i^n + \sqrt{n} \Tilde{x}^n_i) \,\le\,  \Tilde\kappa_2\, n+  
\Bigl(\sum_{j \in \cJ(i)} \mu_{ij}^n \Bigr) \sqrt{n}  \Tilde{x}^n_i \,.
\end{aligned}
\end{equation*}

We next calculate an estimate for $F_{n,i}^{(1)}$ in \cref{PT8C}.
First observe that 
\begin{equation} \label{PT8E} 
\sum_{j\in\cJ(i)}\mu_{ij}^n z_{ij} \,=\,
\sum_{j\in\cJ(i)}\mu_{ij}^n \widetilde{N}_{ij}^n +
\sum_{j\in\cJ(i)}\mu_{ij}^n \bigl(z_{ij}-\widetilde{N}_{ij}^n\bigr)\,.
\end{equation}
We distinguish two cases.

\smallskip\noindent\emph{Case~A.}
Suppose that $x_i\le\widetilde{N}_i^n$.
In this case we have $z_{ij}-\widetilde{N}_{ij}^n\le0$ and
$\Tilde{x}^n_i\le 0$.
Thus we obtain
\begin{equation*}
- \Tilde{x}^n_i\sum_{j\in\cJ(i)}\mu_{ij}^n \bigl(z_{ij}-\widetilde{N}_{ij}^n\bigr)
\,\le\,
- \Tilde{x}^n_i \biggl(\min_{j\in\cJ(i)}\,\mu_{ij}^n\biggr)
\bigl(x_i-\widetilde{N}_i^n) \,=\,
-  \biggl(\min_{j\in\cJ(i)}\,\mu_{ij}^n\biggr)\, \sqrt{n}\abs{\Tilde{x}^n_i}^2\,.
\end{equation*}
Therefore, by \cref{EL5A,PT8E},  we have
\begin{equation*}
F_{n,i}^{(1)}(x) \,\le\,
- \Bigl(\vartheta_p\,p_i\,\sqrt{n} + \sorder\bigl(\sqrt n\bigr)\Bigr)\,\Tilde{x}^n_i
-\sqrt{n}\,\biggl(\min_{j\in\cJ(i)}\,\mu_{ij}^n\biggr)\,\abs{\Tilde{x}^n_i}^{2}\,.
\end{equation*}

\smallskip\noindent\emph{Case~B.}
Suppose that $x_i > \widetilde{N}_i^n$.
In this case, $z_{ij}-\widetilde{N}_{ij}^n\ge0$ and
$\Tilde{x}^n_i\ge 0$.
By  \cref{EL5A,PT8C,PT8E}, we then immediately have that
\begin{equation*}
F_{n,i}^{(1)}(x) \,\le\, \Tilde{x}^n_i
\biggl(\lambda_i^n - \sum_{j \in \cJ(i)} \mu_{ij}^n \widetilde{N}^n_{ij} \biggr)
 \,\le\, - \Bigl(\vartheta_p\,p_i\,\sqrt{n} + \sorder\bigl(\sqrt n\bigr)\Bigr)\,\Tilde{x}^n_i\,.
\end{equation*}
From Cases A--B, we obtain 
\begin{equation*}
\begin{aligned}
F_{n,i}^{(1)}(x)  &\,\le\, 
-\Bigl( \vartheta_p\,p_i\sqrt{n} + \sorder\bigl(\sqrt n\bigr) \Bigr)\,
\Tilde{x}^n_i\,\Ind_{\{\Tilde{x}^n_i>0\}} \\
&\mspace{100mu}
- \sqrt{n}\,\Biggl(\biggl(\vartheta_p\,p_i
+ \frac{\sorder\bigl(\sqrt n\bigr)}{\sqrt n}\biggr)\,\Tilde{x}^n_i
+\biggl(\min_{j\in\cJ(i)}\,\mu_{ij}^n\biggr)\,
\abs{\Tilde{x}^n_i}^{2}\Biggr)\,\Ind_{\{\Tilde{x}^n_i\le0\}}\,.
\end{aligned}
\end{equation*}
Using these estimates, we deduce that
for $\varepsilon>0$ small enough and for all $n\ge n_0$,
there exist positive constants 
$C_k$, $k=0,1$, satisfying
\begin{equation*}
\widetilde{\Lg}_n^z\,\widetilde\Lyap_{\epsilon}(x) \,\le\,  C_0
-C_1\,\widetilde\Lyap_{\epsilon}(x)\qquad\forall\,x\in\ZZ^I_+ \,.
\end{equation*}
Exponential ergodicity follows from this drift inequality.
This completes the proof.
\Halmos\endproof

\begin{remark}\label{R9}
We remark that the results in \cref{T8} can be extended 
for networks with renewal arrivals and exponential service times in the same way
as in  
\citet[Section 3.2]{AHP18}. In particular, we include the 
age process $S^n_i(t)$ of each class-$i$ customers into 
the state descriptor so that $(X^n,S^n)$ is a Markov 
process. We use a Lyapunov function as defined in 
 \citet[Eq.~(3.8)]{AHP18} together with the function 
$\Lyap_{\epsilon}(x)$ in \cref{D6.2}. We can then derive the 
associated Foster-Lyapunov equation by combining the 
calculations in \cref{T8} and those of 
\citet[Theorem~3.1]{AHP18}  related to the age processes. 
The same applies to the transience result for the diffusion-scaled processes 
in \cref{T6}.
We leave the details for the reader to verify.
\end{remark}

\section{Concluding Remarks}
In this paper we have introduced an important 
parameter for multiclass multi-pool networks of any tree 
topology, which plays the same critical role as the safety
staffing parameter in the square-root staffing of 
single-class many-server queues  in the Halfin--Whitt 
regime \citep{H-W-81,w92}.  
Our results show that  the SWSS being positive is 
necessary and sufficient for stabilizability for networks 
with renewal arrivals and exponential service times. We 
conjecture that it is also the necessary and sufficient 
condition in the non-Markovian case (networks with non-exponential service times).
This would require a Markovian description of the system dynamics using
measure-valued processes
(e.g., \citet{Kaspi-11,Kaspi-13,Agh-Ram-19,Agh-Ram-20}).
This is an interesting
open problem for future work.

\section*{Acknowledgment}
We thank the reviewers for their helpful comments that have improved the exposition of the paper. In particular, the authors wish to thank the anonymous reviewer that suggested an alternative approach to prove \cref{T3} which also helped in proving \cref{L2}.

This work was supported in part by an Army Research Office grant W911NF-17-1-0019,
in part by NSF grants DMS-1715210, CMMI-1635410, and DMS/CMMI-1715875,
and in part by the Office of Naval Research through grant N00014-16-1-2956,
and approved for public release under DCN\# 43-6911-20.

\bibliographystyle{informs2014}

\bibliography{MCMP}

\end{document}